\documentclass[11pt]{article}
\usepackage{fullpage}
\usepackage[english]{babel}
\usepackage{amsmath,amssymb,amsfonts}
\newtheorem{theorem}{Theorem}[section]
\newtheorem{lemma}[theorem]{Lemma}
\newtheorem{proposition}[theorem]{Proposition}
\newtheorem{definition}{Definition}[section]

\newtheorem{hypothesis}[theorem]{Hypothesis}
\newtheorem{assumption}[theorem]{Assumption}
\newtheorem{remark}[theorem]{Remark}
\newtheorem{corollary}[theorem]{Corollary}

\def\Dim{\noindent\hbox{{\bf Proof.}$\;\; $}}          
\def\sqr#1#2{{\vcenter{\vbox{\hrule height .#2pt
     \hbox{\vrule width .#2pt height#1pt \kern#1pt \vrule
     width .#2pt} \hrule height .#2pt}}}}
\def\square{\mathchoice\sqr54\sqr54\sqr{4.1}3\sqr{3.5}3}

\makeatletter
\@addtoreset{equation}{section}

\makeatother

\def\qed{{\hfill\hbox{\enspace${ \square}$}} \smallskip}
\def\sqr#1#2{{\vcenter{\vbox{\hrule height .#2pt \hbox{\vrule
 width .#2pt height#1pt \kern#1pt \vrule
width .#2pt} \hrule height .#2pt}}}}
\def\square{\mathchoice\sqr54\sqr54\sqr{4.1}3\sqr{3.5}3}

\def\ds{\begin{displaystyle}}
\def\eds{\end{displaystyle}}
\def\dis{\displaystyle }
\def\<{\langle }
\def\>{\rangle }

\def\R{\mathbb R}

\def\E{\mathbb E}
\def\P{\mathbb P}

\def\calf{{\cal F}}
\def\calg{{\cal G}}

\author{Marco Fuhrman
\\
Dipartimento di Matematica,
Politecnico di Milano\\
piazza Leonardo da Vinci 32, 20133 Milano, Italy\\
e-mail: marco.fuhrman@polimi.it
\\ \\
Ying Hu\\
IRMAR, Universit\'e Rennes 1\\ Campus de Beaulieu, 35042 Rennes
Cedex, France\\ e-mail: ying.hu@univ-rennes1.fr
\\ \\
Gianmario Tessitore
\\
Dipartimento di Matematica e Applicazioni, Universit\`a di Milano-Bicocca\\
Via Cozzi 53, 20125 Milano, Italy\\
e-mail: gianmario.tessitore@unimib.it }

\title{Stochastic Maximum Principle for Optimal Control of
Partial Differential Equations Driven by White Noise}
\date{}
\author{Marco Fuhrman
\\
Dipartimento di Matematica,
Universit\`a di Milano\\
via Saldini 50, 20133 Milano, Italy\\
e-mail: marco.fuhrman@unimi.it
\\ \\
Ying Hu\\
IRMAR, Universit\'e Rennes 1\\ Campus de Beaulieu, 35042 Rennes
Cedex, France\\ e-mail: ying.hu@univ-rennes1.fr
\\ \\
Gianmario Tessitore
\\
Dipartimento di Matematica e Applicazioni, Universit\`a di Milano-Bicocca\\
Via Cozzi 53, 20125 Milano, Italy\\
e-mail: gianmario.tessitore@unimib.it }

\begin{document}

\maketitle

\begin{abstract}

We prove  a stochastic maximum principle of
Pontryagin's type for the optimal control of a stochastic partial differential equation
driven by  white noise in the case when
 the set of control actions is convex. Particular attention is paid
 to well-posedness of the adjoint backward
 stochastic differential equation and the regularity
 properties of its solution with values in infinite-dimensional
 spaces.

\end{abstract}

\section{Introduction}

In this paper we consider an infinite-dimensional
stochastic optimal control problem for a system evolving in a
Hilbert space $H$ and
described by
a state equation of the form
\begin{equation}\label{eq:state}
 dX_t=AX_t\,dt+F(X_t,u_t)\,dt+G(X_t,u_t)\,d W_t,\quad X_0=x\in H,
\end{equation}
where
  $A$ is the infinitesimal generator of a strongly continuous
  semigroup $e^{tA}$ of linear operators,   $W$ a cylindrical
Wiener process in $H$, $F$ and $G$ are suitable
drift and diffusion coefficients, with values in $H$
and $\mathcal{L}(H)$ respectively, depending on
 a control process $u$ taking   values in a set
 $\mathcal {U}$ contained in another Hilbert space $U$.
The cost function is
\begin{equation}\label{eq:cost}
  J(x,u)=\E\int_0^TL(X_t,u_t)\,dt+\E \Phi (X_T),
\end{equation}
for suitable real-valued functions $L$, $\Phi$
(more precise assumptions will be given later).

Our goal in this paper is to give a necessary condition for existence of an optimal control.
This
condition, called stochastic maximum principle (SMP) in the sense
of Pontryagin,
was extensively studied in the finite dimensional case, especially after
the seminal paper by Peng \cite{Pe} which gives a very  general
form of the SMP. There is at present a great interest
in generalizations to the infinite dimensional case, that
was started in \cite{Be}, with particular emphasis on the
application to optimal control of stochastic
partial differential equations (SPDEs): see
 \cite{DuMe}, \cite{DuMe2},
\cite{FuHuTe}, \cite{Gu}, \cite{HuPe2}, \cite{LuZh}, \cite{OkSuZh}, \cite{TaLi}, \cite{Zh}.
The main limitation of the present state of the art
is perhaps that known results only deal with the case when the Wiener process
driving the equation is finite-dimensional or has a trace-class covariance
operator whereas, to our
best knowledge, there is no result for a cylindrical Wiener process $W$.

In this paper our purpose is to establish the SMP for
an evolution equation driven by a cylindrical Wiener process,
in a form suitable for direct application to controlled
SPDEs
driven by white noise.
On the other hand, we will  suppose in this paper that the control domain
$\mathcal {U}$
is convex, which allows
us to apply convex perturbation arguments instead of spike variation
in the deduction of the necessary optimality condition.
Our main result is as follows: under suitable conditions,
in particular differentiability conditions
on $F,G,L,\Phi$,  any optimal control $\bar u$ and
the corresponding trajectory $\bar X$ must satisfy the
SMP inequality
$$ \< \nabla_u [F(\bar{X}_t,\bar{u}_t)] (v-\bar{u}_t), p_t\> +
 \<\nabla_uL(\bar{X}_t,\bar{u}_t),v-\bar{u}_t\>
 + {\rm Tr} \left[q_t^*\left(\nabla_{u}[G(\bar{X}_t,\bar{u}_t)](v- \bar{u}_t)\right) \right]
 \geq 0,\quad v\in \mathcal{U}.
$$
Here $\nabla$ denotes the gradient operator,
$(p,q)$ is a pair of stochastic processes
 taking values, respectively,
  in $H$ and  in the space of Hilbert-Schmidt operators
on $H$,
which we characterize as the unique solution to the
so-called adjoint equation, a linear backward stochastic
differential equation (BSDE) of the form
\begin{equation}\label{lbsdesmp-intro}
    \left\{\begin{array}{lll}
    -dp_t&=&\Big[A^*p_t + \sum_{i=1}^\infty C_i^*(t)q_te_i
    +\nabla_x L({\bar X}_t,{\bar u}_t)\Big]\ dt-\sum_{i=1}^\infty q_te_i\,d\beta_t^i,
    \\
    p_T&=&\nabla_x\Phi({\bar X}_T),
\end{array}\right.
\end{equation}
where $e_i$ is an orthonormal basis of $H$,
$\beta^i_t=\<W_t,e_i\>$ are independent scalar Brownian
motions, and $C_i$ are processes in $\mathcal{L}(H)$ defined
as
$C_i (t)  = \nabla_x [G(\bar{X}_t,\bar{u}_t) e_i  ]  $.
There are some techincal difficulties that we have to face.
 The first one is the occurrence of the first
series in \eqref{lbsdesmp-intro}, which needs to be interpreted
in a suitable way and for which convergence holds in a weak sense in general;
we call this class of BSDEs `singular'.
The second one is the occurrence of the trace operator
in the SMP: for this term to be well defined we need to prove
that the process $q$ satisfies additional regularity results,
in particular it takes values in the space of trace-class
operators on $H$, for which the trace is meaningful.

This paper is organized as follows:  in the next section
we introduce our control problem and formulate the main assumptions;
then we prove that they are verified in the reference case
of the nonlinear controlled heat equation perturbed by noise.
In Section 3 we apply the convex perturbation
argument, we deduce an equation for the first-order variation process
and, under some more assumptions, we prove the SMP taking for granted some results
on singular BSDEs. In Section 4
we make a careful study of a general class of linear singular
BSDEs, proving in particular the trace-class regularity
 and useful duality relations. In the Appendices
we conclude with some reminders and some auxiliary results.

\section{Formulation of the Optimal Control Problem}

Given two real separable Hilbert spaces $X$ and $Y$ by $\mathcal{S}_2(X,Y)$ we denote
 the Hilbert space of Hilbert-Schmidt operators from $X$ to $ Y$ and by $\mathcal{S}_1(X,Y)$ we denote the Banach space of trace class operators (these are instances of the
Schatten-von Neumann classes of operators: some reminders are collected in
  Appendix \ref{appendix:trace-class}). We write $\mathcal{S}_1(X)$,
  $\mathcal{S}_2(X)$ instead of $\mathcal{S}_1(X,X)$, $\mathcal{S}_2(X,X)$.

Given two Banach spaces $D$ and $E$, we say that a mapping $f: D\rightarrow E$ is of class $\mathcal{G}^1(D,E)$ if it is G\^{a}teaux differentiable and $\nabla f: D\rightarrow \mathcal{L}(D,E)$ is strongly continuous (that is, continuous with respect to the strong operator topology).

Finally $H$ is a separable Hilbert space and $(W_t)_{t\geq 0}$ is an
 $H$-valued cylindrical Wiener process. We denote  $(\mathcal{F}_t)_{t\geq 0}$
 the corresponding Brownian filtration, completed in the usual way,
 that   verifies the usual conditions.

By $L^p_{\mathcal{P}}([0,T]\times \Omega,H)$ we denote the Banach space of $H$-valued progressively measurable processes $X$ with $\mathbb{E}\int_0^T |X_s|^pds < \infty$ and by $L^p_{\mathcal{P}}(\Omega,{C}([0,T],H))$ the subspace  of $H$-valued progressively measurable processes $X$ with continuous trajectories satisfying $\mathbb{E}\sup_{s\in[0,T]}|X_s|^p < \infty$.
Here and below we use the symbol $|\cdot|$ to denote a norm when  the corresponding
space is clear from the context, otherwise we use a subscript.

\begin{assumption}\label{assumptions}
\begin{enumerate}
   \item[(i)] $A$ is the generator of a strongly continuous semigroup $e^{tA}$, $t\geq 0$,
   of bounded linear operators in $H$.

   \item[(ii)] $\mathcal{U}$ is a
   convex subset of a separable Banach space $U$.

  \item[(iii)] $F:H\times U\to H$ is a map of class $\mathcal{G}^1(H\times U,H)$.
Moreover, denoting its gradient on $H\times U$ as $(\nabla_x F, \nabla_u F)$,
there exist constants $L,C\geq 0$ such that
  $$|\nabla_x F (x,u)|_{\mathcal{L}(H)}+|\nabla_u F (x,u)|_{\mathcal{L}(U,H)}\leq L,
  \quad |F(0,u)|\leq C,\qquad  x\in H,\; u\in \mathcal{U}.
  $$
  \item[(iv)] $G: H\times U\to L(H)$ satisfies $e^{sA} G(x,u)\in \mathcal{S}_2(H)$  for all $s>0$, $x\in H$, $u\in U$, and   the map $(x,u)\rightarrow e^{sA} G(x,u) $ is of class
      $\mathcal{G}^1(H\times U,\mathcal{S}_2(H))$.
   Moreover there exist constants $L,C\geq 0$ and  $\alpha\in [0,1/2)$ such that
   \begin{eqnarray*}
     |\nabla_x [e^{sA} G(x,u)]|_{\mathcal{L}(H, \mathcal{S}_2(H))}
     + |\nabla_u[ e^{sA} G(x,u)]|_{\mathcal{L}(U, \mathcal{S}_2(H))} &\leq &  L s^{-\alpha}, \\
     |e^{sA}G(0,u)|_{ \mathcal{S}_2(H)} &\leq & C\; s^{-\alpha},
   \end{eqnarray*}
    for  all $x\in H$, $u\in \mathcal{U}$.

\item[(v)]
$L: H\times U \rightarrow \mathbb{R}$ and  $\Phi: H \rightarrow \mathbb{R}$ are bounded, Lipschitz and of class $\mathcal{G}^1$.
\item[(vi)] There exists an orthonormal basis $(e_i)_{i\in \mathbb{N}} \in H$ such  that, for all $i\in \mathbb{N} $ and all $u\in U$, the map $ x \rightarrow G(x,u) e_i $ is of class
    $\mathcal{G}^1(H, H)$. Moreover there exists a constant $L\ge 0$ such that
    $$|\nabla_x[G(x,u) e_i]|_{\mathcal{L}(H)}\leq L,
    \qquad i\in \mathbb{N},\;  x\in H,\; u\in U.
    $$

\item[(vii)] For all $x,y \in H$ the map $ u\rightarrow G(x,u)y$ is
of class $\calg^1(U,H)$ and
there exists a constant $C\ge 0$ such that
$${ |\nabla_{u}[G(x,u)y]v|_H\leq C\, |y|_H |v|_U}, \qquad  x,y \in  H,\; u,v \in U.
$$
\end{enumerate}
\end{assumption}

   Any progressively measurable $\mathcal{U}$-valued process $u$
will be called
an admissible control.

Under the above assumptions, for every admissible control, the state equation (formulated in mild sense):
\begin{equation}\label{eq:state-equation} X_t=e^{tA}x+\int_0^t e^{(t-s)A} F({X}_s,{u}_s)ds+\int_0^t {e^{(t-s)A} G(X_s,u_s)}dW_s\end{equation}
admits a unique solution $X\in L^p_{\mathcal{P}}(\Omega,C([0,T],H))$
for every $p\geq 1$, see \cite{DaZa}.
Moreover the cost associated to the initial datum $x$ and control $u$ is the
well defined real number
$$J(x,u)=\mathbb{E}\int_0^T L(X_t,u_t)dt+\mathbb{E}\Phi(X_T).
$$

\begin{remark}\label{commentassumpt}
\emph{
\begin{enumerate}
\item
Combining points $(iv)$ and $(vii)$ of Assumption \ref{assumptions}
it is easy to check that
for every $s>0$, $x,h\in H$, $u\in U$ and any direction $v\in U$
we have
\begin{equation}\label{derivGu}
    (\nabla_u[e^{sA}G(x,u)]v)\,h= e^{sA}\nabla_u[G(x,u)h]v.
\end{equation}
Similarly, $\nabla_x[e^{sA}G(x,u)]\,e_i= e^{sA}\nabla_x[G(x,u) e_i]$.
\item
For further use we need to introduce a suitable approximation of the derivative
operator $(\nabla_u[e^{sA}G(x,u)]v)$.
Let us denote $\Pi_n$ the orthogonal projection in
$\mathcal{L}(H)$ onto the linear span of $e_1,\ldots,e_n$. Next,
for fixed $x\in H$, $u,v\in U$, let us
define operators in $\mathcal{L}(H)$  setting,
for every  $y\in H$,
$$
\Gamma(x,u,v) y=\nabla_{u}[G(x,u)y]v, \;
\Gamma^n(x,u,v) y=\nabla_{u}[G(x,u)(\Pi_ny)]v
=\sum_{i=1}^n\nabla_{u}[G(x,u)e_i]v\, \left<y,e_i\right>_H.
$$
Note that by $(vi)$  we have
$$
|\Gamma(x,u,v)|_{\mathcal{L}(H)}\le C\,|v|,\quad
|\Gamma^n(x,u,v)|_{\mathcal{L}(H)}\le C\,|v|\sqrt {n},\quad
\Gamma^n(x,u,v)y\to \Gamma(x,u,v)y
$$
in $H$ as $n\to\infty$.  Moreover, recalling \eqref{derivGu}, we have for $s>0$
\begin{equation}\label{approxderiv}
|e^{sA}\Gamma^n(x,u,v)- (\nabla_u[e^{sA}G(x,u)]v)|^2_{\mathcal{S}_2(H)}
=\sum_{i>n}| e^{sA} \nabla_u[G(x,u)e_i]v|_H^2
\downarrow 0,
\end{equation}
since the series $\sum_{i=1}^\infty| e^{sA} \nabla_u[G(x,u)e_i]v|_H^2=
|\nabla_u[e^{sA}G(x,u)]v|^2_{\mathcal{S}_2(H)}$ is convergent.
Finally, given a trace class operator $Q\in\mathcal{S}_1(H)$ it is easy to check that
\begin{equation}\label{approxtraccia}
|\hbox{Tr}(Q\Gamma^n)|\le C\, |v|_U\, |Q|_{\mathcal{S}_1(H)},
\qquad
\hbox{Tr}(Q\Gamma^n)\to \hbox{Tr}(Q\Gamma).
\end{equation}
\end{enumerate}
}
\end{remark}

\subsection{The Reference Example}

Consider the following controlled  stochastic heat equation in $[0,1]$:,
$$
\left\{
\begin{array}{lll}
dX_t(\xi)&=&\dis \frac{\partial^2}{\partial x^2}X_t(\xi)\,dt +b(\xi,X_t(\xi),u_t(\xi))\,dt+
\sigma(\xi,X_t(\xi),u_t(\xi))d {\mathcal W}(t,\xi), \\ \\
X_t(0)&=&X_t(1)=0, \quad t\in [0,T],\\
\\
X_0(\xi)&=&x(\xi), \quad \xi\in[0,1],
\end{array}
\right.
$$
where $b$ and $\sigma:[0,1]\times\R\times \R\to \R$ are given
Borel measurable functions. We assume that
 $b(t,\cdot, \cdot)$ and   $\sigma(t,\cdot,\cdot)$
 are of class $C^1$,  Lipschitz  uniformly with respect to $t$, and
  that $b(\cdot, 0,0)$ and   $\sigma(\cdot, 0, 0)$ are bounded.
In the above equation  $({\mathcal W}(t,\xi) )$, $t\ge 0$, $\xi\in [0,1]$ is a {space time white noise} and by $(\calf_t)_{t\ge 0}$
we denote its natural (completed) filtration.
The set of admissible control actions $\mathcal{U}$ is a convex subset of $U:=L^{2}([0,1])$
and we assume that $\mathcal{U}\subset L^{\infty}([0,1])$.
A control $u$ is a (progressive) process with values in $\mathcal{U}$.
We also introduce the { cost functional}:
$$
J(u)=\E \int_0^T\int_0^1 l(\xi,X_t(\xi),u_t(\xi))\,d\xi\,dt
+\E \int_0^1 h(\xi,X_T(\xi))\,d\xi,
$$
where
$l:[0,1]\times\R\times \R\to \R$,
$h:[0,1]\times\R\to \R$ are given bounded, Borel measurable  functions.
We assume that, for a.e. $\xi\in [0,1]$, $b(\xi, \cdot,\cdot)$ and
$h(\xi,\cdot)$ are of class $C^1$ with bounded derivatives (uniformly with respect to $\xi$).

To reformulate the problem in our general framework we have set $H=L^2([0,1])$ and consider an  $H$-valued
{ cylindrical Wiener process}
$ ( W_t )_{t\geq 0}$.

$A$ is the realization of the second derivative operator in $H$
with Dirichlet boundary conditions. So $\mathcal{D}(A)=H^2([0,1])\cap H^1_0([0,1])$ and $A\phi=\phi''$ for all $\phi \in \mathcal{D}(A)$.

Finally for $x,y\in L^2([0,1])$, $u,v \in L^{\infty}([0,1])$,
$${F(x,u)(\xi)=b(\xi,x(\xi),u(\xi))
,\qquad  [G(x,u)y](\xi)=\sigma(\xi,x(\xi),u(\xi))y(\xi),}$$
$$  L(x,u)=\int_{0}^1l(\xi,x(\xi),u(\xi))d\xi,\quad \Phi(\xi)=\int_{0}^1 h(\xi, X(\xi)) d\xi.
$$
The state equation written in abstract form is
$$dX_t=AX_tdt+F(X_t,u_t)dt+G(X_t,u_t)d W_t,\quad X_0=x,
$$
where $x\in H$ and the solution will evolve in $H$. Moreover
the cost becomes
 $$J(x,u)=\E\int_0^TL(X_s,u_s)ds+\E \Phi (X_T).$$
 It is well known (see \cite{DaZa2}) that Assumptions \ref{assumptions} (\textit{i})-(\textit{v}) are satisfied.
 Concerning Assumption \ref{assumptions} $(vi)$ we notice that in our concrete case for all $y\in L^2([0,1])$,
$$[\nabla_x (G(x,u) e_i  )y](\xi)= \frac{\partial \sigma}{\partial x} (\xi,x(\xi),u(\xi))e_i(\xi)y(\xi).$$
So it is enough to choose an orthonormal basis in $L^2([0,1])$ with $\sup_i \sup_\xi
|e_i(\xi)|<\infty$, for instance a trigonometrical  basis.

Finally concerning Assumption \ref{assumptions} $(vii)$ we notice that
for all $y\in  L^2([0,1])$ and all $v\in L^{\infty}([0,1])$:
$$[\nabla_u ((G(x,u)y)]v(\xi)= \frac{\partial\sigma}{\partial u} (\xi,x(\xi),u(\xi))y(\xi)v(\xi),$$
and $|[\nabla_u ((G(x,u)y)]v(\xi)|_{L^2([0,1])}\leq L_{\sigma} |y|_{L^2([0,1])} |v|_{L^\infty([0,1])}$ where
$L_{\sigma}$ is the Lipschitz constant of
$\sigma$ with respect to $u$.

\section{The Stochastic Maximum Principle}

\subsection{First Variation}
In this section we perturb a given admissible control, that eventually will be the optimal one, and  compute the corresponding expansion of the cost.

In the following  $u$ and $u'$ are two admissible controls and we assume that $\delta u := u'-u$ is bounded in $U$.
Moreover we denote by $X^{\epsilon}$ the state (e.g. the solution to equation \eqref{eq:state-equation}  corresponding to the control $u^{\epsilon}:=(1-\epsilon)u + \epsilon u'=
u+\epsilon \delta u$).

Regular dependence on parameters of the mild solution to forward SDEs gives us the first order expansion of the state:
\begin{theorem}
The map $\epsilon \rightarrow X^{\epsilon}$ is of class $\mathcal{G}^1$ from $[0,1)$ to  $L^p_{\mathcal{P}}(\Omega,C([0,T],H))$ and its derivative in $\epsilon=0$ is given by the unique mild solution $Y\in L^p_{\mathcal{P}}(\Omega,C([0,T],H))$ of the linear equation
\begin{equation}\label{eq:first-variation-Y}
\left\{
\begin{array}{lll}
dY_t&=&
\Big[AY_t +\nabla_x F(X_t,u_t)Y_t + \nabla_uF(X_t,u_t)\delta u_t
\Big]\,dt
\\&&
+\Big[\nabla_x G(X_t,u_t)Y_t + \nabla_u G(X_t,u_t)\delta u_t\Big]\,dW_t,
\\
Y_0&=&0.
\end{array}
\right.
\end{equation}
Explicitly, we have $\P$-a.s.
\begin{equation}\label{firstvariationmild}
    \begin{array}{lll}
Y_t&=&\dis \int_0^te^{(t-s)A}
\nabla_x F(X_s,u_s)Y_s \,ds + \int_0^te^{(t-s)A} \nabla_uF(X_s,u_s)\delta u_s
\,ds
\\&&
\dis +\int_0^t \nabla_x [e^{(t-s)A}  G(X_s,u_s)]Y_s \,dW_s
+\int_0^t
 \nabla_u [e^{(t-s)A}G(X_s,u_s)]\delta u_s\,dW_s,
 \qquad t\in [0,T].
\end{array}
\end{equation}
\end{theorem}
\Dim The proof will follow arguments similar to the ones exposed in \cite{FT1}, see the proof of Proposition 3.3. We limit ourselves to proving the claim in the case $F=0$, the general case being a straightforward extension.

Consider the mapping $\Phi$ from $L^p_{\mathcal{P}}(\Omega,C([0,T],H))\times [0,1)$ to  $L^p_{\mathcal{P}}(\Omega,C([0,T],H))$ given by
$$\Phi(\Xi,\epsilon)_t=e^{tA}x+\int_0^t {e^{(t-l)A} G(\Xi_l,u_l+ \epsilon \delta u_l)}dW_l.$$
Clearly $X^{\epsilon}$ is a solution to the state equation \eqref{eq:state-equation} with $u$ replaced by $u^{\epsilon}$ if and only if it is a fixed point of $\Phi(\cdot,\epsilon)$.

In \cite{FT1} it is shown that, if $\beta >0$ is large enough then $\Phi(\cdot,\epsilon)$ is a contraction, uniformly with respect to $\epsilon$, in  $L^p_{\mathcal{P}}(\Omega,C([0,T],H))$ endowed with the equivalent norm $(\mathbb{E}\sup_{t\in [0,T]} e^{\beta t}|\Xi_t|^p)^{1/p}$. Moreover $\Phi(\cdot,\epsilon)$ is of class $\mathcal{G}^1$ from  $L^p_{\mathcal{P}}(\Omega,C([0,T],H))$ to $L^p_{\mathcal{P}}(\Omega,C([0,T],H))$ with derivative, in the direction $N$, given by:
$$(\nabla_{\Xi} \Phi(\Xi,\epsilon)N)_t= \int_0^t {e^{(t-l)A} \nabla_{x}G(\Xi_l,u_l+\epsilon \delta u_l)}N_ldW_l.$$
Concerning the dependence on $\epsilon$ we have
 \begin{align*}
I^{h}_t:&=\frac{\Phi(\Xi, \epsilon + h)_t- \Phi(\Xi, \epsilon)_t}{h} -\int_0^t \nabla_u \left( e^{(t-l)A}G(\Xi_l,u_l)\delta u_l\right)dW_l \\
&  =\int_0^t\left\{\int_0^1 \left[ \nabla_u\left( e^{(t-l)A}G(\Xi_l,u^{\epsilon}_l+\zeta h \delta u_l)\right)\delta u_l- \nabla_u\left( e^{(t-l)A}G(\Xi_l, u^{\epsilon}_l)\right)\delta u_l\right]d\zeta\right\} dW_l.
\end{align*}
By the factorization method (see the proof of Proposition 3.2 in \cite{FT1}) we get for $1/p< \gamma < 1/2-\alpha$:
$$|I^h|_{L^p_{\mathcal{P}}(\Omega,C([0,T],H))}^p\leq c_p \mathbb{E}\int_0^T|V^h_l|^p dl,$$
where
$$V^h_l= \int_0^l\!\! (l-\sigma)^{-\gamma}\left\{\int_0^1 \left[ \nabla_u\left( e^{(l-\sigma)A}G(\Xi_\sigma,u^{\epsilon}_\sigma+\zeta h \delta u_\sigma)\right)\delta u_\sigma- \nabla_u\left( e^{(l-\sigma)A}G(\Xi_\sigma, u^{\epsilon}_\sigma)\right)\delta u_\sigma\right]d\zeta\right\} dW_\sigma.$$
By the Burkholder-Davis-Gundy inequality
\begin{align*}\mathbb{E}|V^h_l|^p \leq c_p
\mathbb{E} \left\{\int_0^l(l-\sigma)^{-2\alpha}\left[\int_0^1 \left| \nabla_u\left( e^{(l-\sigma)A}G(\Xi_\sigma,u^{\epsilon}_\sigma+\zeta h \delta u_\sigma)\right)\delta u_\sigma \right.\right.\right.\qquad\quad \\
\left.\left.\left. -\nabla_u\left( e^{(l-\sigma)A}G(\Xi_\sigma,u^{\epsilon}_\sigma)\right)\delta u_\sigma\right|^2_{\mathcal{S}_2(K,H)} d\zeta \right] d\sigma \right\}^{p/2}.
\end{align*}
Thus by Hypothesis \ref{assumptions}-$(iv)$:
$$\mathbb{E}|V^h_l|^p \leq c_p \left\{|\delta u|^2_{L^{\infty}(\Omega\times [0,T],U)}\int_0^l (l-\sigma)^{-2(\gamma+\alpha)} d\sigma\right\}^{p/2}\leq c_p l^{p/2- p(\gamma+\alpha)}|\delta u|^p_{L^{\infty}(\Omega\times [0,T],U)}.$$
The continuity of $\nabla_u\left( e^{(s-\sigma)A}G(\Xi, \cdot)\right)$ and Dominated Convergence Theorem yield:
$\mathbb{E}|Y^h_s|^p \rightarrow 0$ and consequently that $|I^h|_{L^p_{\mathcal{P}}(\Omega,C([0,T],H))}\rightarrow 0$. We can therefore conclude that $\Phi$ is differentiable with respect to $\epsilon$ as well, with
$$\nabla_{\epsilon} \Phi (\Xi,\epsilon)=\int_0^t \nabla_u \left( e^{(t-l)A}G(\Xi_l,u_l)\right)\delta u_ldW_l.$$
The continuity of $\nabla_\Xi \Phi (\Xi,\epsilon)$  and  $\nabla_\epsilon \Phi (\Xi,\epsilon)$ with respect to $\Xi $ and $\epsilon$ can be proved in a similar way.

Summing up, $\Phi$ is a mapping of class $\mathcal{G}^1$ on $L^p_{\mathcal{P}}(\Omega,C([0,T],H))\times [0,1)$. The parameter depending contraction principle (see \cite{FT1} Proposition 2.4) yields that the map $[0,1)\ni\epsilon \rightarrow X^{\epsilon}\in L^p_{\mathcal{P}}(\Omega,C([0,T],H)) $ that relates the parameter to the fixed point is of class $\mathcal{G}^1$ (in this case just differentiable with continuous derivative). Moreover its derivative satisfies:
$$\nabla_{\epsilon}X^{\epsilon}=\nabla_{\Xi}\Phi(X^\epsilon, \epsilon)\nabla_{\epsilon}X^{\epsilon}+\nabla_{\epsilon}\Phi(X^\epsilon, \epsilon).$$
Plugging in the above relation the expressions for
$\nabla_{\Xi}\Phi$ and $\nabla_{\epsilon}\Phi$ we get,
$$\nabla_{\epsilon}X^{\epsilon}_t
=\int_0^t e^{(t-l)A} \nabla_{x}G(X^\epsilon_l,u_l+\epsilon \delta u_l)
\nabla_{\epsilon} X^{\epsilon}_l\, dW_l
+\int_0^t e^{(t-l)A} \nabla_{u}G(X^\epsilon_l,u_l+\epsilon \delta u_l)\delta u_l\,dW_l,$$
and the claim follows letting $\epsilon=0$ and denoting $Y_t$ the limit. \qed

As a consequence of the above result we have the following expansion of the cost:
\begin{proposition}\label{expansioncost} With the above notation, we have:
$$J(x,u^{\epsilon})= J(x,{u})+{\epsilon I(\delta u)}+o(\epsilon),$$
where
$$ {I(\delta u)}=\E\int_0^T\Big[\<\nabla_x L({X}_t,{u}_t),Y_t\>
+\<\nabla_uL({X}_t,{u}_t),\delta u_t\>\Big]dt+
\E\<\nabla_x\Phi({X}_T),Y_T\>.$$
\end{proposition}
\Dim By the above theorem, if
$R^{\epsilon}_t:=\epsilon^{-1}[X^{\epsilon}_t-X_t -\epsilon Y_t]$ then $|R^{\epsilon}|_{L^p_{\mathcal{P}}(\Omega,C([0,T],H))}\rightarrow 0$.
Moreover
\begin{align*}
J(x,u^{\epsilon})-J(x,u)&=\mathbb{E}\int_0^T \left[L(X_t+\epsilon (Y_t+R^\epsilon_t), u_t+ \epsilon \delta u_t)- L(X_t,u_t)\right] dt\\& \quad +\mathbb{E}\left[\Phi(X_T+\epsilon (Y_T+R^\epsilon_T))- \Phi(X_T)\right]\\
&=\epsilon \mathbb{E}\int_0^T \int_0^1\nabla_x L(X_t+\lambda\epsilon( Y_t+R^\epsilon_t), u_t+ \lambda\epsilon \delta u_t)(Y_t+R^\epsilon_t)
\, d\lambda\, dt \\
& \quad + \epsilon \mathbb{E}\int_0^T \int_0^1\nabla_u L(X_t+\lambda\epsilon( Y_t+R^\epsilon_t), u_t+ \lambda\epsilon \delta u_t) \delta u_t \, d\lambda\,dt\\
&\quad +\epsilon \mathbb{E}\int_0^1\nabla_x \Phi(X_T+\lambda \epsilon( Y_T+R^\epsilon_T))( Y_T+R^\epsilon_T)
\, d\lambda,
\end{align*}
taking into account the continuity and boundedness of $\nabla_x L$, $\nabla_u L$ and $\nabla_x \Phi$,
applying the dominated convergence theorem it is then easy to show that
\begin{align*}
&\mathbb{E}\int_0^T \int_0^1\nabla_x L(X_t+\lambda\epsilon( Y_t+R^\epsilon_t), u_t+ \lambda\epsilon \delta u_t)Y_t
\, d\lambda\, dt \rightarrow \mathbb{E}\int_0^T \nabla_x L(X_t, u_t)Y_t
\, dt, \\
& \mathbb{E}\int_0^T \int_0^1\nabla_x L(X_t+\lambda\epsilon( Y_t+R^\epsilon_t), u_t+ \lambda\epsilon \delta u_t)R^{\epsilon}_t
\, d\lambda\, dt \rightarrow 0, \\ &
 \mathbb{E}\int_0^T \int_0^1\nabla_u L(X_t+\lambda\epsilon( Y_t+R^\epsilon_t), u_t+ \lambda\epsilon \delta u_t) \delta u_t \, d\lambda\,dt \rightarrow  \mathbb{E}\int_0^T \nabla_u L(X_t, u_t) \delta u_t \,dt,\\
 &\mathbb{E}\int_0^1\nabla_x \Phi(X_T+\lambda \epsilon( Y_T+R^\epsilon_T))Y_T \rightarrow \mathbb{E}\nabla_x \Phi(X_T) Y_T, \\ &\mathbb{E}\int_0^1\nabla_x \Phi(X_T+\lambda \epsilon( Y_T+R^\epsilon_T))R^{\epsilon}_T \rightarrow 0,
\end{align*}
and the proof is completed. \qed







\subsection{Stochastic Maximum Principle}

We will prove the stochastic maximum principle under
the following additional assumption:
\begin{assumption}\label{semigrouphs}
We have
$e^{tA}\in \mathcal{S}_2(H)$ for all $ t> 0$ and there exist $c>0$ and
$\alpha < 1/2 $ such that
 $$|e^{tA}|_{\mathcal{S}_2(H)} \leq c t^{-\alpha},\qquad  t\in  (0,T].
$$
\end{assumption}

In order to state the stochastic maximum principle
we assume that
 an optimal control $\bar{u}$   exists and we denote $\bar{X}$ is the corresponding state.
Next we need to introduce
the dual process $(p_t,q_t)$, with values in $H\times
\mathcal{S}_2(H)$.
To this end
we fix an orthonormal basis $\{e_i\}_{i\in \mathbb{N}}$ in $H$ such that point $(vi)$ in Assumption \ref{assumptions} holds and we
define $C_i (t) h= \nabla_x [G(\bar{X}_t,\bar{u}_t) e_i  ]h $ for $h\in H$.  Then we have
\begin{equation}\label{cilimit}
        |C_i(t)|_{\mathcal{L}(H)}\leq L.
\end{equation}
Recalling Remark \ref{commentassumpt}-1 and
taking into account Hypothesis \ref{assumptions}-$(iv)$ we also obtain
\begin{equation}\label{cihs}
\sum_{i=1}^{\infty} |e^{tA}C_i(s) h|^2 =
\sum_{i=1}^{\infty} |\nabla_x[e^{tA}G  (\bar{X}_s,\bar{u}_s) e_i  ]h|^2
=|\nabla_x[e^{tA}G  (\bar{X}_s,\bar{u}_s)  ]h|^2_{\mathcal{S}_2(H)}
\leq L^2 t^{-2\alpha}|h|_H^2,
\end{equation}
for all $ t> 0,$ $s\geq 0$, $h\in H$.

Next we introduce the adjoint equation for the unknown process
$(p,q)$,
written formally
\begin{equation}\label{lbsdesmp}
    \left\{\begin{array}{lll}
    -dp_t&=&\Big[A^*p_t + \sum_{i=1}^\infty C_i^*(t)q_te_i
    +\nabla_x L({\bar X}_t,{\bar u}_t)\Big]\ dt-\sum_{i=1}^\infty q_te_i\,d\beta_t^i,
    \\
    p_T&=&\nabla_x\Phi({\bar X}_T),
\end{array}\right.
\end{equation}
where $\beta^i_t=\< e_i, W_t\>$, $i=1,2...$ is a family of independent Brownian motions.
 The precise notion of (mild) solution to this
 equation is as follows: we say that a pair
  $(p,q)$ with  $p\in L^2_{{\cal P}}(\Omega\times [0,T],H)$, $q\in L^2_{{\cal P}}(\Omega\times [0,T],\mathcal{S}_2(H))$ is a mild solution to equation \eqref{lbsdesmp} if,
  for any $t\in [0,T]$, we have
\begin{equation}\label{eq:mildsmp}
\begin{array}{lll}
p_t&=&\dis
e^{(T-t)A^*}\nabla_x\Phi({\bar X}_T)
+ \sum_{i=1}^\infty \int_t^T e^{(s-t)A^*} C_i^*(s)Q_se_i    \,ds
\\&&\dis
+\int_t^T e^{(s-t)A^*}\nabla_x L({\bar X}_s,{\bar u}_s)\, ds-\sum_{i=1}^\infty\int_t^T e^{(s-t)A^*}Q_s e_i\,d\beta_s^i,\qquad \mathbb{P}-a.s.,
\end{array}
\end{equation}
where, for fixed $t$, the series
$ \sum_{i=1}^\infty \int_t^T e^{(s-t)A^*} C_i^*(s)Q_se_i    \,ds$
is required to converge  weakly in the space
$ L^2(\Omega,{\cal F}_T,\mathbb P,H)$
[in Proposition \ref{prop:strong} below we will also give additional conditions
that guarantee that the series converges in a stronger sense].
We refer the reader to the next section, in particular to
Definition \ref{def:solution_BSDE} and Remark \ref{defweakconv},
for a more precise discussion of this notion of solution.
There we will also prove the following result.

\begin{proposition}\label{existuniqsmp}
Under Assumptions \ref{assumptions}
  and \ref{semigrouphs} there exists a unique solution $(p,q)$ to
  equation \eqref{lbsdesmp}.
Moreover, $q_t\in \mathcal{S}_1(H)$   $ d\mathbb{P}\otimes dt$-a.s. and
$$\mathbb E\int_0^T (T-t)^{2\alpha}|q_t|^2_{\mathcal{S}_1(H)}dt<\infty .
$$
\end{proposition}

\noindent {\bf Proof.} This follows from
Theorem \ref{Th:existence} and
Proposition \ref{prop:trace-class-reg}.
  Hypothesis \ref{Hyp:C_i}, needed for these statements
 to hold, is verified due to  Assumption \ref{semigrouphs}
 and inequalities \eqref{cilimit} and \eqref{cihs}.
\qed

The final ingredient in the proof of the stochastic maximum principle
is a duality relation involving the first variation process $Y$ solution
to equation \eqref{eq:first-variation-Y} with $u=\bar{u}$ and $X=\bar X$.
In the present notation this equation (to be understood in the mild
form \eqref{firstvariationmild})
can be written formally as
\begin{equation}\label{eq:first-variation-Y-modif}
\left\{
\begin{array}{lll}
dY_t&=&
\Big[AY_t +\nabla_x F(\bar X_t,\bar u_t)Y_t + \nabla_uF(\bar X_t,\bar u_t)\delta u_t
\Big]\,dt
\\&&\dis
+\sum_{i=1}^{\infty} C_i(t)Y_t \,d\beta^i_t +
\sum_{i=1}^{\infty} \nabla_u [G(\bar X_t,\bar u_t)e_i]\delta u_t
\,d\beta^i_t,
\\
Y_0&=&0.
\end{array}
\right.
\end{equation}
 In the next section we will prove the following duality relation.

\begin{proposition}\label{dualitysmp}
With the previous assumptions and notations,
suppose that
$\rho:[0,T]\times \Omega\to H$ and
$\Gamma: [0,T]\times \Omega \rightarrow \mathcal{S}_2(H)$
are progressively measurable and bounded
and $\mathcal{Y}$ denotes the unique mild solution to the equation
\begin{equation}\label{eq:first-variation-Y-modif-2}
\left\{
\begin{array}{lll}
d\mathcal{Y}_t&=&
\Big[A\mathcal{Y}_t +\nabla_x F(\bar X_t,\bar u_t)\mathcal{Y}_t + \rho_t
\Big]\,dt
+\sum_{i=1}^{\infty} C_i(t)\mathcal{Y}_t \,d\beta^i_t +
\sum_{i=1}^{\infty} \Gamma_t e_i
\,d\beta^i_t,
\\
\mathcal{Y}_0&=&0.
\end{array}
\right.
\end{equation}
Then
\begin{equation}\label{eq:dualityMsmp}
\mathbb E\int_0^T\langle p_t,\rho_t\rangle dt+\mathbb E\int_0^T\langle q_t,\Gamma_t\rangle_{\mathcal{S}_2(H)}dt
=\mathbb E \langle \nabla_x\Phi({\bar X}_T),\mathcal{Y}_T\rangle
+\mathbb E \int_0^T\langle \nabla_x L({\bar X}_t,{\bar u}_t), \mathcal{Y}_t\rangle dt.
\end{equation}
\end{proposition}

\noindent\textbf{Proof.}
This is exactly formula
\eqref{eq:dualityM} of Corollary  \ref{Cor:construction-solution},
where we put $s=0$, $x=0$,  $\gamma=0$,
$\eta=\nabla_x\Phi({\bar X}_T)$, $f_t=\nabla_x L({\bar X}_t,{\bar u}_t)$
and note that in this case
$\tilde{\mathcal{Y}}^{\infty,M}$ coincides with the solution $\mathcal{Y}$
to \eqref{eq:first-variation-Y-modif-2}.
\qed

Now we are ready to state and prove the main result of this paper.

\begin{theorem}\label{SMP}
Suppose that  Assumptions \ref{assumptions}
  and \ref{semigrouphs} hold and that  an  optimal pair $(\bar{u}, \bar{X})$
  exists. Then for every
 $  v\in \mathcal{U}$ we have,  $ d\mathbb{P}\otimes dt$-a.s.,
$$ \< \nabla_u [F(\bar{X}_t,\bar{u}_t)] (v-\bar{u}_t), p_t\> +
 \<\nabla_uL(\bar{X}_t,\bar{u}_t),v-\bar{u}_t\>+ {\rm Tr} \left[q_t^*\left(\nabla_{u}[G(\bar{X}_t,\bar{u}_t)](v- \bar{u}_t)\right) \right]\geq 0,
$$
where $(p,q)$ is the  unique mild solution to
  equation \eqref{lbsdesmp}.
\end{theorem}

\noindent {\bf Proof.} \emph{Step 1:} we prove the duality formula
\begin{equation}\label{eq:dualityMsmpexact}
\begin{array}{l}\dis
\mathbb E\int_0^T\langle p_t,\nabla_uF(\bar X_t,\bar u_t)\delta u_t\rangle dt
+\mathbb E\int_0^T{\rm Tr}
\left[q_t^*\left(\nabla_{u}[G(\bar{X}_t,\bar{u}_t)]\delta u_t\right) \right]\,dt
\\\dis
=\mathbb E \langle \nabla_x\Phi({\bar X}_T),{Y}_T\rangle
+\mathbb E \int_0^T\langle \nabla_x L({\bar X}_t,{\bar u}_t), {Y}_t\rangle dt,
\end{array}
\end{equation}
where $Y$ is the first variation process solution to \eqref{eq:first-variation-Y-modif}.

We define
$$\rho_t= \nabla_uF(\bar X_t,\bar u_t)\delta u_t, \qquad
 \Gamma_th= \nabla_u [G(\bar X_t,\bar u_t)h]\delta u_t,
 $$
for every $h\in H$. Since we take $\delta u$ to be a bounded process,
it follows that $\rho $ is also bounded, by Assumption
\ref{assumptions}-$(iii)$.
Heuristically, we note that with this choice
the equations \eqref{eq:first-variation-Y-modif} and \eqref{eq:first-variation-Y-modif-2}
coincide, so that $Y=\mathcal{Y}$ and \eqref{eq:dualityMsmpexact}
coincides with \eqref{eq:dualityMsmp}.
However such argument is not correct, as we can not directly apply
Proposition \ref{dualitysmp} above, since
$\Gamma$ is not a bounded process with values
in $\mathcal{S}_2(H)$, so we have to revert to an approximation
procedure.

Let us denote $\Pi_n$ the orthogonal projection in
$\mathcal{L}(H)$ onto the linear span of $e_1,\ldots,e_n$ and define
$$
 \Gamma^n_th= \nabla_u [G(\bar X_t,\bar u_t)(\Pi_nh)]\delta u_t
 =\sum_{i=1}^n\nabla_{u}[G(\bar X_t,\bar u_t)e_i]\delta u_t\, \left<h,e_i\right>_H,
\qquad h\in H.
$$
Each $\Gamma^n$ is a bounded process
in $\mathcal{L}(H)$ (by Assumption \ref{assumptions}-$(vii)$
and since $\delta u$ is bounded) and since
it  has finite rank it is also bounded in $\mathcal{S}_2(H)$.
 Let $\mathcal{Y}^n$ be the unique mild solution to the equation
\begin{equation}\label{eq:first-variation-Y-modif-3}
\left\{
\begin{array}{lll}
d\mathcal{Y}^n_t&=&
\Big[A\mathcal{Y}^n_t +\nabla_x F(\bar X_t,\bar u_t)\mathcal{Y}^n_t + \rho_t
\Big]\,dt
+\sum_{i=1}^{\infty} C_i(t)\mathcal{Y}^n_t \,d\beta^i_t +
\sum_{i=1}^{\infty} \Gamma^n_t e_i
\,d\beta^i_t,
\\
\mathcal{Y}_0^n&=&0.
\end{array}
\right.
\end{equation}
We can now apply
Proposition \ref{dualitysmp} and obtain the duality relation
\begin{equation}\label{eq:dualityMsmpapprox}
\mathbb E\int_0^T\langle p_t,\rho_t\rangle dt+\mathbb E\int_0^T\langle q_t,\Gamma^n_t\rangle_{\mathcal{S}_2(H)}dt
=\mathbb E \langle \nabla_x\Phi({\bar X}_T),\mathcal{Y}^n_T\rangle
+\mathbb E \int_0^T\langle \nabla_x L({\bar X}_t,{\bar u}_t), \mathcal{Y}^n_t\rangle dt.
\end{equation}
Now we let $n\to\infty$. It is convenient at this point to recall the notation
introduced in Remark \ref{commentassumpt}-2, namely the operators
$\Gamma(x,u,v)$ and $\Gamma^n(x,u,v)$: indeed we have
$$
\Gamma_t= \Gamma({\bar X}_t,{\bar u}_t,\delta u_t),\qquad
\Gamma^n_t= \Gamma^n({\bar X}_t,{\bar u}_t,\delta u_t),
$$
and it follows from
\eqref{approxderiv} that
for $s>0$
\begin{equation}\label{approxnablau}
    |e^{sA}\Gamma^n_t- (\nabla_u[e^{sA}G({\bar X}_t,{\bar u}_t)]
\delta u_t)|^2_{\mathcal{S}_2(H)}=
|e^{sA}\Gamma^n({\bar X}_t,{\bar u}_t,\delta u_t)- (\nabla_u[e^{sA}G({\bar X}_t,{\bar u}_t)]\delta u_t)|^2_{\mathcal{S}_2(H)}
\downarrow 0,
\end{equation}
$d\P\otimes dt$-a.s.,
and from \eqref{approxtraccia} that
\begin{equation}\label{approxtrace}
|\hbox{Tr}(q_t^*\Gamma^n_t)|\le C\, |\delta u_t|_U\, |q_t|_{\mathcal{S}_1(H)}
\le c\,  |q_t|_{\mathcal{S}_1(H)},
\quad
\hbox{Tr}(q_t^*\Gamma^n_t)\to \hbox{Tr}(q_t^*\Gamma_t),
\qquad
d\P\otimes dt-a.s.,
\end{equation}
since we know that
$q_t\in\mathcal{S}_1(H)$ $d\P\otimes dt$-a.s.
Writing down the mild form of the equations for $Y$ and $\mathcal{Y}^n$
and substracting we obtain
(compare \eqref{firstvariationmild}),
   $$ \begin{array}{lll}
\mathcal{Y}^n_t-Y_t&=&\dis \int_0^te^{(t-s)A}
\nabla_x F(X_s,u_s)(\mathcal{Y}^n_s-Y_s) \,ds
+\int_0^t \nabla_x [e^{(t-s)A}  G(X_s,u_s)](\mathcal{Y}^n_s-Y_s) \,dW_s
\\&&
\dis
+\int_0^t\bigg(e^{(t-s)A}\Gamma^n_t-
 \nabla_u [e^{(t-s)A}G(X_s,u_s)]\delta u_s\bigg)\,dW_s.
\end{array}
$$
It follows from \eqref{approxnablau} that the last integral tends to zero
in $L^2(\Omega, \calf,\P,H)$. By standard estimates it also follows
that $\sup_{t\in [0,T]}\E |\mathcal{Y}^n_t-Y_t|_H^2\to 0$,
so that we can pass to the limit in the right-hand side of \eqref{eq:dualityMsmpapprox}.
Next we note that
$$
\mathbb E\int_0^T\langle q_t,\Gamma^n_t\rangle_{\mathcal{S}_2(H)}dt
=
\mathbb E\int_0^T\hbox{Tr}(q_t^*\Gamma^n_t)dt
\to \mathbb E\int_0^T\hbox{Tr}(q_t^*\Gamma_t)dt
$$
by dominated convergence, as  it follows
from \eqref{approxtrace} and the fact that
$$
\mathbb E\int_0^T
   |q_t|_{\mathcal{S}_1(H)}\, dt \le
  \bigg( \int_0^T(T-t)^{-2\alpha}dt\bigg)^{1/2}
\bigg(\mathbb E\int_0^T
(T-t)^{2\alpha}   |q_t|^2_{\mathcal{S}_1(H)}\, dt\bigg)^{1/2} <\infty
$$
by Proposition \ref{existuniqsmp} and the assumption that $\alpha <1/2$.
Passing to the limit in \eqref{eq:dualityMsmpapprox} we
finish the proof of Step 1.

 \emph{Step 2:} conclusion.
It follows from
Proposition \ref{expansioncost}
and the duality formula of Step 1 that
\begin{eqnarray*}
J(x,u^{\epsilon})-J(x,\bar{u})  &  =&
\epsilon\,\mathbb{E}\int_0^T\< \nabla_u F(\bar{X}_t,\bar{u}_t)\delta u_t, p_t\>\, dt+
\epsilon \mathbb{E}\int_0^T \<\nabla_uL(\bar{X}_t,\bar{u}_t),\delta u_t\>\,dt\\
 &&   +  \epsilon\mathbb{E}\int_0^T{ \hbox{Tr}
 \left[q_t^*\left(\nabla_{u}[G(\bar{X}_t,\bar{u}_t)]\delta u_t\right) \right]} \,dt + o(\epsilon)
\end{eqnarray*}
Since $\bar u$ is optimal, we have $J(x,u^{\epsilon})-J(x,\bar{u})\ge 0$
and the proof can be concluded by standard arguments
based on  localization and the  Lebesgue differentiation theorem, see, e.g., \cite{Pe, YoZh}.
\qed

\section{Singular  Infinite-dimensional BSDEs}

The main purpose of this section is to give a complete
proof of
Propositions \ref{existuniqsmp} and \ref{dualitysmp} that
were used in an essential way to prove the stochastic maximum principle.
They both refer to properties of the dual (backward)
equation \eqref{lbsdesmp}. To simplify the notation we will
present our results  in the case when $F=0$, the general case
being essentially the same.
On the other hand, we will address a class of backward equations
which are otherwise more general, namely of the form
\begin{equation}\label{lbsde}
    \left\{\begin{array}{lll}
    -dP_t&=&[A^*P_t + \sum_{i=1}^\infty C_i^*(t)Q_te_idt+f_t]\, dt-\sum_{i=1}^\infty Q_te_i\,d\beta_t^i,
    \\
    P_T&=&\eta,
\end{array}\right.
\end{equation}
where $\beta^i_t=\< e_i, W_t\>$, $i=1,2...$ is a family of independent Brownian motions,
  $\eta\in L^2(\Omega,{\cal F}_T,\mathbb P,H)$, $f\in L^2_{\cal P}(\Omega\times [0,T],H)$.
The unknown process is the pair denoted $(P,Q)$ and takes values
in $H\times \mathcal{S}_2(H)$.
We will work under the following assumptions, which are assumed to
hold throughout this section.

\medskip

\begin{hypothesis}\label{Hyp:C_i} $ $
\begin{enumerate}

  \item  $e^{tA}$, $t\ge0$, is a strongly continuous semigroup of
  bounded linear operators in $H$. Moreover,
  $ e^{tA}\in {\mathcal{S}_2(H)}$ for all $ t> 0$ and there exist
  constants $c>0$ and $\alpha \in [0, 1/2) $ such that
  $|e^{tA}|_{\mathcal{S}_2(H)} \leq c t^{-\alpha}$ for all $ t\in ( 0,T]$.

\item The processes $C_i $ are strongly progressively
measurable with values in $\mathcal{L}(H)$. Moreover
we have $ |C_i(t)|_{\mathcal{L}(H)}\leq c,\quad \mathbb{P}-\hbox{a.s.}$ for all $t\in [0,T]$
and $i\in \mathbb{N}$.

  \item $\sum_{i=1}^{\infty} |e^{tA}C_i(s) h|^2 \leq c t^{-2\alpha}|h|_H^2$ for
  all $ t\in ( 0,T],$ $s\geq 0$, $h\in H$.

\end{enumerate}
\end{hypothesis}

\medskip

Strongly progressively measurable processes means that
they are progressively measurable with values in
the space $\mathcal{L}(H)$ endowed with
  the Borel sets of the strong operator topology.

We notice that the sum $\sum_{i=1}^\infty C_i^*(t)Q_te_i$
in equation \ref{lbsde}
is  not convergent  in general, even when $Q_t$ is Hilbert-Schmidt. Moreover the semigroup does not seem to be directly helpful since
under Hypothesis \ref{Hyp:C_i} it is not clear whether the sum
$\sum_{i=1}^\infty e^{sA^*}C_i^*(t)Q_te_i$ converges or not.
 The main result of this paper is the proof of well-posedness of such class of linear BSDEs
 which are driven by white noise and involve a `singular' infinite sum.

 We give the following notion of (mild) solution.

 \begin{definition}\label{def:solution_BSDE} We say that a pair of processes $(P,Q)$ with  $P\in L^2_{{\cal P}}(\Omega\times [0,T],H)$, $Q\in L^2_{{\cal P}}(\Omega\times [0,T],\mathcal{S}_2(H))$ is a mild solution to equation \eqref{lbsde} if the following holds:

  \begin{enumerate}

\item Denoting   $S^M(s) : =\sum_{i=1}^M (T-s)^{\alpha}  C_i^*(s)Q_se_i$, $s\in [0,T]$ then the sequence $(S^M)$ converges weakly in $L^{2}_{\mathcal{P}}(\Omega \times [0,T],H)$.

\item For any $s\in [0,T]$,
\begin{equation}\label{eq:mild}P_s=e^{(T-s)A^*}\eta+ \sum_{i=1}^\infty \int_s^T e^{(l-s)A^*} C_i^*(l)Q_le_idl+\int_s^T e^{(l-s)A^*}f_ldl-\sum_{i=1}^\infty\int_s^T e^{(l-s)A^*}Q_l e_id\beta_l^i,\, \mathbb{P}-a.s.\end{equation}
\end{enumerate}
 \end{definition}

\begin{remark}\label{defweakconv}{\em
Notice that, for any fixed $s$, the map $g\mapsto \int_s^T (T-l)^{-\alpha}e^{(l-s)A^*}g_l dl$ is a bounded linear functional from $L^2_{\mathcal{P}}([0,T]\times \Omega,H)$ to $ L^2(\Omega,{\cal F}_T,\mathbb P,H)$, hence weakly continuous. So  if condition 1 above holds  then,
for all fixed $s\in [0,T]$, the sum
$$ \sum_{i=1}^M \int_s^T e^{(l-s)A^*} C_i^*(l)Q_ldl=\int_s^T (T-l)^{-\alpha}e^{(l-s)A^*} S^M(l)dl$$
converges, weakly in $ L^2(\Omega,{\cal F}_T,\mathbb P,H)$,  to a limit
 that we denote $\sum_{i=1}^\infty \int_s^T  e^{(l-s)A^*} C_i^*(l)Q_le_i dl$  and
 that appears in \eqref{eq:mild}.
} \end{remark}

\subsection{Linear Forward SDEs}
We will study equation \eqref{lbsde} exploiting duality arguments. To
this end we start  by collecting  precise estimates on the solutions of a suitable family
of linear
forward SDEs.
Namely, given any starting time $s\in [0,T]$, we consider the equation on the time interval $[s,T]$:
\begin{equation}\label{lfsde}
    \left\{\begin{array}{lll}
    d\mathcal{Y}_t&=&A\mathcal{Y}_t\,dt + \sum_{i=1}^\infty C_i(t)\mathcal{Y}_t\,d\beta_t^i+\sum_{i=1}^\infty C_i(t)\gamma_t\,d\beta_t^i  + \sum_{i=1}^\infty \Gamma_t e_i\, d\beta_t^i+\rho_tdt,
    \\
    \mathcal{Y}_s&=&x,
\end{array}\right.
\end{equation}
together with the approximating equations, for $N,M\in \mathbb{N}\cup \{\infty\}$:
\begin{equation}\label{lfsden}
    \left\{\begin{array}{lll}
    d\tilde{\mathcal{Y}}_t^{N,M}&=&A\tilde{\mathcal{Y}}^{N,M}_t\,dt + \sum_{i=1}^N C_i(t)\tilde{\mathcal{Y}}^{N,M}_t\,d\beta_t^i+\sum_{i=1}^M C_i(t)\gamma_t\,d\beta_t^i  + \sum_{i=1}^\infty \Gamma_t e_i\, d\beta_t^i+\rho_tdt,
    \\
    \tilde{\mathcal{Y}}^{N,M}_s&=&x.
\end{array}\right.
\end{equation}
In the above equation we always assume that $x:\Omega \rightarrow H$ is bounded and $\mathcal{F}_s$ measurable, $\rho,\gamma:[s,T]\times \Omega\to H$
are progressively measurable and bounded, $\Gamma: [s,T]\times \Omega \rightarrow \mathcal{S}_2(H)$  is also  progressively measurable and bounded.

For further use we note that $    \tilde{\mathcal{Y}}^{N,M}$ clearly does
not depend on $M$ when $\gamma=0$.

We start from a standard estimate on this SDE. Its proof   coincides with the one of Proposition 3.2 in \cite{FT1} and will be omitted.

\begin{theorem}\label{thm:lfsde}
For all $p\in [2,\infty)$ and
 $N,M\in \mathbb{N}\cup \{\infty\}$, in the space  $L^p_{\mathcal{P}}(\Omega,C([s,T],H))$
there exists a unique solution $\mathcal{Y}$ to equation (\ref{lfsde}) and   a unique solution $\tilde{\mathcal{Y}}^{N,M} $ to equation (\ref{lfsden}); note that $\tilde{\mathcal{Y}}^{\infty,\infty}={\mathcal{Y}}$. Moreover, if  $p$ is large enough, the following estimate holds:
$$\mathbb E(\sup_{t\in [s,T]}|\mathcal{Y}_t|^p)\le c_p\left(1+\mathbb E |x|^p+|\Gamma|_{L^\infty([s,T]\times \Omega,\mathcal{L}(H))}^p+|\gamma|_{L_{\cal P}^\infty([s,T]\times \Omega,H)}^p+\mathbb E\left(\int_s^T|\rho_t|dt\right)^p\right).$$
We also have
\begin{equation}\label{convforward}
    \tilde{{\mathcal{Y}}}^{N,M}\rightarrow \tilde{{\mathcal{Y}}}^{N,\infty},
\quad \tilde{{\mathcal{Y}}}^{N,M}\rightarrow \tilde{{\mathcal{Y}}}^{\infty,M},
\quad
\tilde{{\mathcal{Y}}}^{M,M}\rightarrow {{\mathcal{Y}}},
\quad
\tilde{{\mathcal{Y}}}^{N,\infty}\rightarrow {{\mathcal{Y}}},
\quad
\tilde{{\mathcal{Y}}}^{\infty,M}\rightarrow {{\mathcal{Y}}}
\end{equation}
in the norm of $L^p_{\mathcal{P}}(\Omega,{C}([s,T],H))$, namely
$\mathbb E(\sup_{t\in [s,T]}| \tilde{\mathcal{Y}}_t^{N,M}
-\tilde{\mathcal{Y}}_t^{N,\infty}|^p)\rightarrow 0$
etc.

In addition,
the above  estimate holds for $\tilde{\mathcal{Y}}^{N,M}$,
 $\tilde{{\mathcal{Y}}}^{N,\infty}$,  $\tilde{{\mathcal{Y}}}^{\infty,M}$
 uniformly with respect to $N$ and $M$.

Finally, the first (respectively, the second) convergence result in \eqref{convforward}
holds true uniformly with respect to  $N$ (respectively, to $M$).
\end{theorem}
The next estimate involves the Hilbert-Schmidt norm of $\Gamma$.

\begin{proposition}\label{estimate1}
Under the above assumptions and notations it holds
$$\mathbb \sup_{t\in [s,T]} \mathbb{E}|\mathcal{Y}_t|^2\le c\mathbb E\left[|x|^2+\left(\int_s^T|\rho_t|dt\right)^2+
|\gamma|^2_{L_{\cal P}^\infty([s,T]\times \Omega,H)}
+\int_s^T|\Gamma_t|_{\mathcal{S}_2(H)}^2dt\right].$$
\end{proposition}
\noindent{\bf Proof.} Writing equation \eqref{lfsde} in the mild form, namely
$$\mathcal{Y}_t=e^{(t-s)A}x+\int_s^te^{(t-l)A}\rho_l d l +\sum_{i=1}^{\infty}\int_s^te^{(t-l)A}\left[C_i(l)(\mathcal{Y}_l+\gamma_l)+\Gamma_l e_i\right]d\beta^i_l,
$$
and taking into account Hypothesis \ref{Hyp:C_i} we get
\begin{eqnarray}
\nonumber
  \mathbb{E}|\mathcal{Y}_t|^2 &\le &
  c\mathbb E\left[|x|^2+\left(\int_s^t|\rho_l|dl\right)^2
+\int_s^t (t-l)^{-2\alpha}|\gamma_l|^2dl
+\int_s^t|e^{(t-l)A}\Gamma_l|_{\mathcal{S}_2(H)}^2dl\right] \\
    & & +c\int_s^t(t-l)^{-2\alpha}\mathbb{E}|\mathcal{Y}_l|^2 d l.
    \label{stimafor}
\end{eqnarray}
We obtain
$$\mathbb{E}|\mathcal{Y}_t|^2\le c\mathbb E\left[|x|^2+\left(\int_s^T|\rho_l|dl\right)^2
+|\gamma|_{L^\infty}^2+\int_s^T|\Gamma_l|_{\mathcal{S}_2(H)}^2dl\right] +c\int_s^t(t-l)^{-2\alpha}\mathbb{E}|\mathcal{Y}_l|^2 d l,
$$
and the claim then follows applying a variant of the Gronwall Lemma, see Lemma
7.1.1 in \cite{Henry}.
\qed

The final estimate will be an important tool in the rest of the paper and again  exploits
a special version of the Gronwall Lemma.

\begin{proposition}\label{estimate2} For $s\le t\le T$,
\begin{equation}\label{eq:gronwall-estimate}\mathbb E|\mathcal{Y}_t|^2\le c\mathbb{E}\left[|x|^2+\int_s^t(t-l)^{-2\alpha}|\gamma_l|_H^2dl
+\int_s^t(t-l)^{-2\alpha}|\Gamma_l|_{\mathcal{L}(H)}^2dl+
 \left(\int_s^t|\rho_l|dl\right)^2\right].\end{equation}
\end{proposition}
\noindent{\bf Proof.} We first note  that
 $$|e^{(t-l)A}\Gamma_l |_{\mathcal{S}_2(H)}|\leq |e^{(t-l)A}|_{\mathcal{S}_2(H)}|\Gamma_l|_{\mathcal{L}(H)}\leq c(t-l)^{-\alpha}|\Gamma_l|_{\mathcal{L}(H)}.$$
Consequently, letting
$$u(t):=\mathbb E|\mathcal{Y}_t|^2, \quad
v(t):=c\mathbb E(|\Gamma_t|^2_{L(H)}+|\gamma_t|^2),\quad w(t):=c\mathbb E|x|^2+c\mathbb E\left(\int_s^t|\rho_l|dl\right)^2,  $$
it follows from \eqref{stimafor} that
$$u(t)\le w(t)+\int_s^t(t-l)^{-2\alpha}v(l)dl+c\int_s^t(t-l)^{-2\alpha}u(l)dl= e(t)+c\int_s^t(t-l)^{-2\alpha}u(l)dl,$$
where we set $e(t)=w(t)+\int_s^t(t-l)^{-2\alpha}v(l)dl.$ Using
again the Gronwall Lemma in \cite{Henry} Lemma 7.1.1 we obtain
\begin{equation}\label{pregronwall}
    u(t)\le e(t)+c\int_s^t(t-l)^{-2\alpha}e(l)dl.
\end{equation}
Next we note that
\begin{eqnarray*}
\int_s^t(t-l)^{-2\alpha}e(l)dl
&=&\int_s^tw(l)(t-l)^{-2\alpha}dl+c\int_s^t(t-l)^{-2\alpha}\int_s^l(l-\zeta)^{-2\alpha}v(\zeta)d\zeta\\
&\le & w(t)\int_s^t(t-l)^{-2\alpha}dl
+c\int_s^t(t-l)^{-2\alpha}\int_s^l(l-\zeta)^{-2\alpha}v(\zeta)d\zeta\\
&\le &cw(t)+c\int_s^t v(\zeta)\int_\zeta^t(t-l)^{-2\alpha}(l-\zeta)^{-2\alpha}dl\, d\zeta.
\end{eqnarray*}
But since
\begin{align*}
\int_\zeta^t(t-l)^{-2\alpha}(l-\zeta)^{-2\alpha}dl &=\int_0^{t-\zeta}(t-\zeta-s)^{-2\alpha}s^{-2\alpha}ds
=\int_0^1(1-z)^{-2\alpha}z^{-2\alpha}(t-\zeta)^{-4\alpha+1}dz\\
&= c(t-\zeta)^{-4\alpha+1}\le c(t-\zeta)^{-2\alpha},
\end{align*}
from  \eqref{pregronwall}  we obtain the required conclusion:
$$u(t)\le cw(t)+c\int_s^t (t-l)^{-2\alpha}v(l)dl.
$$
\qed

\begin{corollary} Integrating the inequality \eqref{eq:gronwall-estimate}
with respect to $t$ we get
\begin{equation}\label{eq:gronwall-estimate-int}
 E\int_s^T|\mathcal{Y}_l|^2 dl\le c \mathbb{E}\left[|x|^2+\int_s^T\left(|\Gamma_l|_{L(H)}^2+|\gamma_l|^2\right)dl+\left(\int_s^T|\rho_l|dl\right)^2\right].\end{equation}
\end{corollary}

\subsection{Existence of a Solution to the Singular BSDE}

We will proceed by approximation. Namely, we consider the following BSDE in infinite dimensions where the singular sum in the drift has been truncated:
\begin{equation}\label{lbsdeN}
    \left\{\begin{array}{lll}
    -dP^N_s&=&[A^*P^N_s + \sum_{i=1}^N C_i^*(s)Q^N_te_i+f_s]\,  ds-\sum_{i=1}^\infty Q^N_se_i\,d\beta_s^i,
    \\
    P^N_T&=&\eta.
\end{array}\right.
\end{equation}
We still assume that
 $\eta\in L^2(\Omega,{\cal F}_T,\mathbb P,H)$, $f\in L^2_{\cal P}(\Omega\times [0,T],H)$,
so the above equation fits in the classical theory of Hilbert valued BSDEs, see \cite{HuPe}.
In particular it holds:
\begin{proposition}\label{prop:solutionN}
Assume that
 $\eta\in L^2(\Omega,{\cal F}_T,\mathbb P,H)$, $f\in L^2_{\cal P}(\Omega\times [0,T],H)$.
Then there exists a unique $(P^N,Q^N)$ with
$$P^N\in L^2_{\cal P}(\Omega,C([0,T],H),\quad Q^N\in L^2_{\cal P}(\Omega\times [0,T],\mathcal{S}_2(H))$$
verifying (\ref{lbsdeN}) in the following mild sense:
\begin{equation}
\label{eq:BSDEmildN}
 P^N_s=e^{(T-s)A^*}\eta+\int_s^T e^{(l-s)A^*}[\sum_{i=1}^N C_i^*(l)Q^N_le_i+f_l]\ dl-\sum_{i=1}^\infty \int_s^T e^{(l-s)A^*} Q^N_le_i\,d\beta_l^i.\end{equation}
Moreover choosing as before $x:\Omega \rightarrow H$ bounded and $\mathcal{F}_s$ measurable, $\rho$, $\gamma$ from $[s,T]\times \Omega$ to $H$ progressively measurable and bounded and $\Gamma: [s,T]\times \Omega \rightarrow \mathcal{S}_2(H)$  progressively measurable and bounded it holds:
\begin{equation}\label{eq:dualityNM}
\begin{array}{l}
 \displaystyle \mathbb E\langle P_s^N,x\rangle+\mathbb E\int_s^T\langle P_l^N,\rho_l\rangle dl+\mathbb E\int_s^T\langle Q_l^N,\Gamma_l\rangle_{\mathcal{S}_2(H)}dl+\mathbb E\int_s^T\sum_{i=1}^M\langle Q_l^Ne_i,C_i(l)\gamma_l\rangle dl\\
\displaystyle =\mathbb E \langle\eta,
\tilde{\mathcal{Y}}^{N,M}_T\rangle +\mathbb E\int_s^T\langle f_l,\tilde{\mathcal{Y}}_l^{N,M}\rangle dl,
\end{array}
\end{equation}
where $\tilde{\mathcal{Y}}^{N,M}$ is the solution to equation \eqref{lfsden}.
\end{proposition}

We now define a candidate solution $(P,Q)$ as the weak limit of $(P^N,Q^N)$ in some Hilbert space.

\begin{corollary} \label{Cor:construction-solution} We have
\begin{enumerate}
\item $P^N$ converges weakly to an element $P$ in $L^2_\mathcal{P}([0,T]\times\Omega,H)$;

\item for any $t$, $P_t^N$ converges weakly to an element $\tilde{P}_t$ in $L^2(\Omega, {\cal F}_t,\mathbb P, H)$;

\item $Q^N$ converges weakly to an element $Q$ in
$L^2_\mathcal{P}([0,T]\times \Omega, \mathcal{S}_2(H))$.
\end{enumerate}

Moreover choosing $\eta$, $f$, $x$, $\rho$, $\gamma$, $\Gamma$ as in Proposition     \ref{prop:solutionN}  it holds, for all $M\in \mathbb{N}$:
\begin{equation}\label{eq:dualityM}
\begin{array}{l}
\displaystyle\mathbb E\langle \tilde{P}_s,x\rangle+\mathbb E\int_s^T\langle P_l,\rho_l\rangle dl+\mathbb E\int_s^T\langle Q_l,\Gamma_l\rangle_{\mathcal{S}_2(H)}dl
+\mathbb E\int_s^T\sum_{i=1}^M\langle Q_le_i,C_i(l)\gamma_l\rangle dl
\\
\displaystyle =\mathbb E \langle\eta,\tilde{\mathcal{Y}}^{\infty,M}_T\rangle +\mathbb E \int_s^T\langle f_l,\tilde{\mathcal{Y}}^{\infty,M}_l\rangle dl.
\end{array}
\end{equation}
\end{corollary}

\noindent\textbf{Proof.}
Let us consider the processes  $\mathcal{Y}$ and
$\tilde{\mathcal{Y}}^{N,M}$, defined as solutions to equations
\eqref{lfsde} and \eqref{lfsden}.
In the first part of the proof we take $s=0$, $x=0$
and $\gamma=0$ in these equations and we recall that
the processes $\tilde{\mathcal{Y}}^{N,M}$  do not depend on $M$,
so that  in particular
$\tilde{\mathcal{Y}}^{\infty,M}= \mathcal{Y}$.

By the estimates in Proposition \ref{estimate1}  the maps
$$\mathcal T:(\rho,\Gamma)\mapsto (\tilde{\mathcal{Y}}^{\infty,M}_T,\tilde{\mathcal{Y}}^{\infty,M}),
\qquad \mathcal T_N:(\rho,\Gamma)\mapsto (\tilde{\mathcal{Y}}^{N,M}_T,\tilde{\mathcal{Y}}^{N,M})
$$
can be extended to   bounded linear maps from the space
$ L^2_{{\cal P}}(\Omega\times [0,T],H)\times L^2_{{\cal P}}(\Omega\times [0,T],\mathcal{S}_2(H)) $ to
$ L^2(\Omega,{\cal F}_T,\mathbb P,H)\times L^2_{\mathcal{P}}(\Omega \times [0,T],H)$.
We denote by $\mathcal T^*, \mathcal T^*_N$ their Hilbert space adjoints. Given
arbitrary
$\eta\in L^2(\Omega,{\cal F}_T,\mathbb P,H)$, $f\in L^2_{\cal P}(\Omega\times [0,T],H)$
and setting $(P,Q)=\mathcal T^*(\eta,f)$, we see that
$$
\mathbb E\int_0^T\langle P_l,\rho_l\rangle dl+\mathbb E\int_0^T\langle Q_l,\Gamma_l\rangle_{\mathcal{S}_2(H)}dl=\mathbb E \langle\eta,\mathcal{Y}^{M}_T\rangle +\mathbb E \int_0^T\langle f_l,\mathcal{Y}^M_l\rangle dl,
$$
while \eqref{eq:dualityNM} (with $s=0$, $x=0$, $\gamma=0$)
shows that $(P^N,Q^N)=\mathcal T_N^*(\eta,f)$.
Since $\tilde{{\mathcal{Y}}}^{N,M}\rightarrow \tilde{{\mathcal{Y}}}^{\infty,M}$
as specified in Theorem \ref{thm:lfsde}, it follows easily that
$P^N\to P$ weakly in $L^2_\mathcal{P}([0,T]\times\Omega,H)$
 and $Q^N\to Q$ weakly in $L^2_\mathcal{P}([0,T]\times \Omega, \mathcal{S}_2(H))$.

Now, for arbitrary $s\in [0,T]$ and $x\in L^2(\Omega,{\cal F}_s,\mathbb P,H)$
it follows from  \eqref{eq:dualityNM} (still with $\gamma=0$)
that  $\mathbb E\langle P_s^N,x\rangle$
has a limit as $N\to\infty$, equal to
$$
-\mathbb E\int_s^T\langle P_l^N,\rho_l\rangle dl-\mathbb E\int_s^T\langle Q_l^N,\Gamma_l\rangle_{\mathcal{S}_2(H)}dl
+\mathbb E \langle\eta,
\tilde{\mathcal{Y}}^{N,M}_T\rangle +\mathbb E\int_s^T\langle f_l,\tilde{\mathcal{Y}}_l^{N,M}\rangle dl.
$$
This shows that $P_s^N$ converges weakly  in $L^2(\Omega, {\cal F}_s,\mathbb P, H)$,
and we denote by $\tilde{P}_s$ its limit.

Finally, coming back to general $\gamma$, equality \eqref{eq:dualityM} follows from \eqref{eq:dualityNM} letting $N\rightarrow \infty$. \qed

 We then need to get some  regularity property of the process $(\tilde{P}_s)$ which, for the moment,  we can not identify with $(P_s)$ and has only been defined $\mathbb P$-a.s. for any fixed $s\in[0,T]$.
\begin{proposition}
The map $s\mapsto \tilde{P}_s$ is weakly continuous from $[0,T]$ to $ L^2(\Omega,{\cal F}_T,\mathbb P,H)$.
\end{proposition}
\noindent{\bf Proof.} In order to stress its dependence on the initial time and state, given  $s\in [0,T]$ and $x\in L^2(\Omega,{\cal F}_s,\mathbb P,H)$, denote  by $\mathcal{Y}^{x,s}$ the solution of equation \eqref{lfsde} with $\rho=\Gamma=\gamma=0$, namely:
\begin{equation}
\label{eq:forward-simple}
 d\mathcal{Y}^{\xi,s}_t=A\mathcal{Y}^{\xi,s}_tds+\sum_{i=1}^\infty C_i(s)\mathcal{Y}^{\xi,s}_t d\beta_s^i,\quad \mathcal{Y}^{\xi,s}_s=x.
\end{equation}
We know by \cite{DaZa} that the above equation admits a unique mild solution with
$\sup_{t\in [s,T]} \mathbb{E}|\mathcal{Y}^{x,s}_t|^2 \leq c_{T}(1+\mathbb E|x|^2)$, morever if $x'\in L^2(\Omega,{\cal F}_s,\mathbb P,H)$ then $\sup_{t\in [s,T]} \mathbb{E}|\mathcal{Y}^{x,s}_t-\mathcal{Y}^{x',s}_t|^2 \leq c_{T}(\mathbb E|x-x'|^2)$.

For fixed $x\in L^2(\Omega,{\cal F}_T,\mathbb P,H)$,
let us denote $x_s:=\mathbb E(x|{\cal F}_s)$. Then equation \eqref{eq:dualityM} yields:
$$\mathbb E\langle \tilde{P}_s,x\rangle=\mathbb E\langle \tilde{P}_s,x_s\rangle=\mathbb E\langle \mathcal{Y}_T^{x_s,s},\eta\rangle+\mathbb E\int_s^T\langle \mathcal{Y}^{x_s,s}_l,f_l\rangle dl.
$$
Since
$$\sup_{s\in [0,T]}\sup_{t\in [s,T]}\mathbb E|\mathcal{Y}_t^{s,x_s}|^2\le \sup_{s\in [0,T]} c_T(1+\mathbb E|x_s|^2)\leq c,$$
the weak continuity is proved if we show that,
for all $ t >s$, the map $s\mapsto \mathcal{Y}^{s,x_s}_t$ is continuous
in the norm of $L^2(\Omega,{\cal F},\mathbb P,H)$.

If $s_n\downarrow {s}$ then
$\mathbb{E}|\mathcal{Y}^{s_n,x_{s_n}}_t-\mathcal{Y}_t^{{s},x_{{s}}}|^2\le 2
\mathbb{E}|\mathcal{Y}^{s_n,x_{s_n}}_t-\mathcal{Y}_t^{{s}_n,x_{{s}}}|^2
+2\mathbb{E}|\mathcal{Y}^{s_n,x_s}_t-\mathcal{Y}_t^{{s},x_{{s}}}|^2$.
The first term is controlled by $\mathbb E|x_{s_n}-x_{s}|^2=\mathbb E|\mathbb E(x|\mathcal{F}_{s_n})-\mathbb{E}(x|\mathcal{F}_{s})|^2\rightarrow 0$.
\smallskip

The second term is $\mathbb{E}|\mathcal{Y}^{s_n,x_s}_t-\mathcal{Y}^{s_n,\mathcal{Y}^{s,x_s}_{s_n}}_t|^2$ and is controlled by
 \begin{eqnarray*}
\mathbb E |\mathcal{Y}^{s,x_s}_{s_n}-x_{s}|^2
&\le& 2\mathbb E|e^{(s_n-{s})A} x_{{s}}-x_{{s}}|^2+ 2\mathbb E\left|\sum_{i=1}^\infty \int_{{s}}^{s_n}e^{(s_n-t)A}C_i(t)\mathcal{Y}_{t}^{{s},x_{{s}}}d\beta_t^i\right|^2\\
&\le&2\mathbb E|e^{(s_n-{s})A} x_{{s}}-x_{{s}}|^2+ c \sup_{l\in [s,T]} \mathbb E|\mathcal{Y}_{l}^{s,x_s}|^2 \int_{s}^{s_n}(s_n -t)^{-2\alpha}dt ,
\end{eqnarray*}
which converges to $0$ as $s_n\downarrow 0$.

On the other hand, if $s_n\uparrow s$,
$\mathcal{Y}_t^{s_n,x_{s_n}}-\mathcal{Y}_t^{s,x}
=\mathcal{Y}_t^{s,\mathcal{Y}^{s_n,x_{s_n}}_s}-\mathcal{Y}_t^{s,x}.$
Hence it is enough to prove that
$\mathcal{Y}^{s_n,x_{s_n}}_s\rightarrow x.$ in the norm of $L^2(\Omega,{\cal F},\mathbb P,H)$.
But, proceeding as before, we have
$$\mathbb{E}|\mathcal{Y}_s^{x_{s_n},s_n}- x|^2
\le 2\mathbb{E}|e^{(s-s_n)A}x_{s_n}-x|^2+c\int_{s_n}^s (s-l)^{-2\alpha}\mathbb{E}|\mathcal{Y}^{x_{s_n},s_n}_l|^2 dl.$$
The claim follows
recalling that $\sup_n\sup_{l\in [s_n,s]}\mathbb E|\mathcal{Y}^{x_{s_n},s_n}|^2
\le c_T \sup_n (1+\mathbb E |x_{s_n}|^2)
<\infty$ and that
$|e^{(s-s_n)A}x_{s_n}-x|\leq c|x_{s_n}-x|+|e^{(s-s_n)A}x-x|.$
\qed

\medskip
\noindent
We are now in a position to prove that $\tilde{P}$ and $P$ coincide.
\begin{proposition} $\tilde P$ is a progressively measurable process and     $\tilde P=P$,
 $dt\otimes d\mathbb{P}$-a.s.
\end{proposition}
\noindent{\bf Proof.} We first prove progressive measurability.  Fixed an arbitrary $t\in [0,T]$ we choose a basis $\{\varphi_m\}$ in $L^2(\Omega, {\cal F}_t,\mathbb P,H)$. We have
$\tilde{P}_l=\sum_{m=1}^\infty (\mathbb E\langle \tilde{P}_l,\varphi_m\rangle_H) \varphi_m$, $\forall l\leq t$. Since $\mathbb E\langle \tilde{P}_l,\varphi_m\rangle_H$ is a continuous  function we immediately deduce that $\tilde{P}$ restricted to $[0,t]$ is $\mathcal{B}[0,t] \otimes \mathcal{F}_t$ measurable.

\medskip
\noindent To show that $(\tilde{P})$ and $(P)$ coincide choose $x=\gamma=\Gamma=0$  and an arbitrary bounded progressively measurable process $\rho$ in \eqref{lfsde}. By definition of $\tilde{P}$, for all $t\in [0,T]$, $\mathbb{E} \langle P^N_t,\rho_t\rangle \rightarrow \mathbb{E} \langle \tilde{P}_t,\rho_t\rangle$ so, by dominated convergence theorem (exploiting the measurability of $(\tilde{P})$))
$\int_0^T \mathbb E \langle  P_l^N,\rho_l\rangle dl\rightarrow  \int_0^T \mathbb E \langle \tilde{P}_l,\rho_l\rangle dl.$ But we already know, see Corollary \ref{Cor:construction-solution}, that $\int_0^T \mathbb E \langle  P_l^N,\rho_l\rangle dl\rightarrow  \int_0^T \mathbb E \langle {P}_l,\rho_l\rangle dl$ and the  claim is proved. \qed

\medskip
\noindent In the following when we refer to the process $P$ we will always refer to its version $\tilde{P}$.

\medskip

\noindent We come now to the study of process $Q$.

\begin{lemma}\label{conv-sum_q}
Setting $S^M(s)=\sum_{i=1}^M (T-s)^\alpha C_i^*(s)Q_{s}e_i$,
 $S^M$ converges weakly in $L^2_{\mathcal{P}}(\Omega\times [0,T];H)$. The limit will be denoted by
$\sum_{i=1}^\infty (T-\cdot)^\alpha C_i^*(\cdot)Q_{\cdot}e_i$.\end{lemma}
\textbf{Proof.} Given an arbitrary bounded progressively measurable process $\gamma$ in $H$,
let $\hat{\mathcal{Y}}^{M,\gamma}$ be the mild solution of the equation
\begin{equation}\label{eq:YMgamma}
    \left\{\begin{array}{lll}
    d\hat{\mathcal{Y}}^{M,\gamma}_t&=&A\hat{\mathcal{Y}}^{M,\gamma}_t\,dt + \sum_{i=1}^\infty C_i(t)\hat{\mathcal{Y}}^{M,\gamma}_t\,d\beta_t^i+\sum_{i=1}^M C_i(t)(T-t)^{\alpha}\gamma_t\,d\beta_t^i,
    \\
    \hat{\mathcal{Y}}^{M,\gamma}_0&=&0.
\end{array}\right.
\end{equation}
Similarly, let
$\hat{\mathcal{Y}}^{\gamma}$ be the mild solution of   equation of
\begin{equation}\label{eq:Ygamma}
    \left\{\begin{array}{lll}
    d\hat{\mathcal{Y}}^{\gamma}_t&=&A\hat{\mathcal{Y}}^{\gamma}_t\,dt + \sum_{i=1}^\infty C_i(t)\hat{\mathcal{Y}}^{\gamma}_t\,d\beta_t^i+\sum_{i=1}^\infty C_i(t)(T-t)^{\alpha}\gamma_t\,d\beta_t^i,
    \\
    \hat{\mathcal{Y}}^{\gamma}_0&=&0.
\end{array}\right.
\end{equation}
We see that  equation \eqref{eq:YMgamma} (respectively,  \eqref{eq:Ygamma})
coincide with equation
\eqref{lfsden} (respectively, \eqref{lfsde}) with $s=x=\Gamma=\rho=0$, $N=\infty$ and $\gamma$ replaced by $(T-\cdot)^{\alpha}\gamma$.

The equality \eqref{eq:dualityM} reads:
\begin{equation}\label{converggamma}
      \mathbb{E}\int_0^T\langle S^M(l), {\gamma}_l\rangle dl=
  \mathbb{E}\int_0^T\langle (T-l)^\alpha\sum_{i=1}^M C_i^*(l)Q_le_i,{\gamma}_l\rangle dl
=\mathbb E\langle \hat{\mathcal{Y}}^{M,\gamma}_T,\eta\rangle+\mathbb E\int_0^T \langle \hat{\mathcal{Y}}^{M,\gamma}_l,f_l\rangle dl.
\end{equation}
By the estimate \eqref{eq:gronwall-estimate}, we have  $|  \mathbb{E}\int_0^T\langle S^M(l), {\gamma}_l\rangle dl|\leq c |\gamma|_{L^2_{\mathcal{P}}([0,T]\times \Omega,H)}$. Since the set of bounded
 elements $\gamma$ is dense in $L^2_{\mathcal{P}}([0,T]\times \Omega,H)$
 it follows that
$|   S^M|_{L^2_{\mathcal{P}}([0,T]\times \Omega,H)}\leq c$.
Since, by
Theorem \ref{thm:lfsde},  the right-hand side of \eqref{converggamma}
converges as  $M\to\infty$ (to the limit
$\mathbb E\langle \hat{\mathcal{Y}}^{\gamma}_T,\eta\rangle+\mathbb E\int_0^T \langle \hat{\mathcal{Y}}^{\gamma}_l,f_l\rangle dl$)  when $\gamma$
is bounded, we conclude that
 $S^M$  converges weakly in $L^2_{\mathcal{P}}([0,T]\times \Omega,H)$.
 \qed

\medskip
\noindent If we replace $Q$ by the approximating operators $Q^M$  we obtain the same limit:
\begin{lemma}\label{conv-sum_qM}
$(T-\cdot)^\alpha\sum_{i=1}^M C_i^*(\cdot)Q^M_{\cdot}e_i$ converges weakly to
$(T-\cdot)^\alpha\sum_{i=1}^\infty C_i^*(\cdot)Q_{\cdot}e_i$ in $L^2_{\mathcal{P}}([0,T]\times \Omega,H)$.
\end{lemma}

\noindent{\bf Proof.} The proof of the existence of the weak limit of the sequence
$(T-\cdot)^\alpha\sum_{i=1}^M C_i^*(\cdot)Q^M_{\cdot}e_i$ follows the same argument as in the proof of the above Lemma replacing equation \eqref{eq:YMgamma} by
\begin{equation}\label{eq:YMMgamma}
    \left\{\begin{array}{lll}
    d\hat{\hat{\mathcal{Y}}}^{M,\gamma}_t&=&A\hat{\hat{\mathcal{Y}}}^{M,\gamma}_t\,dt + \sum_{i=1}^M C_i(t)\hat{\hat{\mathcal{Y}}}^{M,\gamma}_t\,d\beta_t^i+\sum_{i=1}^M C_i(t)(T-t)^{\alpha}\gamma_t\,d\beta_t^i,
    \\
   \hat{\hat{\mathcal{Y}}}^{M,\gamma}_0&=&0.
\end{array}\right.
\end{equation}
and replacing the second equality in \eqref{converggamma} by
$$\mathbb \mathbb{E}\int_0^T\langle (T-l)^\alpha\sum_{i=1}^M C_i^*(l)Q^M_le_i,{\gamma}_l\rangle dl=\mathbb E\langle \hat{\hat{\mathcal{Y}}}^{M,\gamma}_T,\eta\rangle+\mathbb E\int_0^T \langle \hat{\hat{\mathcal{Y}}}^{M,\gamma}_l,f_l\rangle dl,$$
which follows from \eqref{eq:dualityNM}.
The proof that the limit equals $(T-\cdot)^\alpha\sum_{i=1}^\infty C_i^*(\cdot)Q_{\cdot}e_i$ comes from the observation that $\mathbb{E}\sup_{s\in[0,T]}|\hat{\hat{\mathcal{Y}}}^{M,\gamma}_s-\hat{\mathcal{Y}}^{\gamma}_s|^p\rightarrow 0$,
which corresponds to the convergence $\tilde{{\mathcal{Y}}}^{M,M}\rightarrow {{\mathcal{Y}}}$
in  Theorem \ref{thm:lfsde}.
\qed

We are now in a position to prove existence of a solution to the singular BSDE \eqref{lbsde}. For the uniqueness part we need the following Lemma on linear BSDEs with unbounded forcing term which proof will be postponed to the Appendix.
\begin{lemma}\label{lemma:exploding}
Assume that $\xi$ is a progressively measurable process in $H$ with $\mathbb E\int_0^T (T-l)^{2\alpha}|\xi_l|^2dl<\infty$. Then for any $N\in\mathbb N$ and any $\eta\in L^2(\Omega,\mathcal{F}_T,\mathbb{P})$
there exists a unique pair of processes $(p,q)$ with $p$    progressively measurable
in $H$, with $s\mapsto p_s$ continuous from $ [0,T]$ to $L^2(\Omega, \mathcal{F}, H)$,
 and $q\in L^2_{{\cal P}}(\Omega\times [0,T],\mathcal{S}_2(H))$ such that:
\begin{equation}\label{eq:exploding}p_s=e^{(T-s)A^*}\eta+\int_s^T e^{(l-s)A^*}\sum_{i=1}^{N}C_i^*(l)q_le_idl+\int_s^T e^{(l-s)A^*}\xi_l dl -\sum_{i=1}^\infty\int_s^T e^{(l-s)A^*}q_le_id\beta_l^i.\end{equation}
Moreover, letting   $\gamma=\rho=0$, $M=\infty$ and
$\Gamma\in L^\infty_{\mathcal{P}}(\Omega\times [0,T],\mathcal{S}_2(H))$  in equation \eqref{lfsden}  the following duality relation holds:
\begin{eqnarray*}
\mathbb E\langle p_s,x\rangle +\mathbb E\int_s^T \langle q_l,\Gamma_l\rangle_{\mathcal{S}_2(H)}ds
&=&\mathbb E\langle \tilde{\mathcal{Y}}^{N}_T,\eta\rangle +\mathbb E\int_s^T \langle (T-l)^{-\alpha}\tilde{\mathcal{Y}}^{N}_l,(T-l)^{\alpha}\xi_l\rangle dl.
\end{eqnarray*}
Notice that $\sup_{l\in [s,T]}\mathbb E|\tilde{\mathcal{Y}}^{N}_l|^2 < +\infty$ and that $\alpha<1/2$, so the last integral is well defined.
\end{lemma}
\noindent{\bf Proof.}
 The proof of this lemma is postponed to the Appendix.\qed
\begin{theorem}\label{Th:existence}
The pair  $(P,Q)$ constructed in Corollary \ref{Cor:construction-solution} is the unique  mild solution  to the singular BSDEs \eqref{lbsde}. \end{theorem}

\noindent\textbf{Proof.} \textit{Existence:} As noticed
in Remark  \ref{defweakconv},   for any fixed $s$, the map $g\mapsto \int_s^T (T-l)^{-\alpha}e^{(l-s)A^*}g_l dl$ is weakly continuous from
$L^2_{\mathcal{P}}([0,T]\times \Omega,H)$ to $ L^2(\Omega,{\cal F}_T,\mathbb P,H)$.

Hence, by Lemma \ref{conv-sum_q}  and Lemma \ref{conv-sum_qM} the sequence
 $ \sum_{i=1}^N \int_s^T e^{(l-s)A^*} C_i^*(l)Q^N_ldl$
converges, weakly in $ L^2(\Omega,{\cal F}_T,\mathbb P,H)$,  to  $\sum_{i=1}^\infty \int_s^T  e^{(l-s)A^*} C_i^*(l)Q_le_i dl$.

Similarly $ \sum_{i=1}^\infty \int_s^T e^{(l-s)A^*} Q^N_le_id\beta_l^i$ converges weakly to $ \sum_{i=1}^\infty \int_s^T e^{(l-s)A^*} Q_le_id\beta_l^i$, since $Q^N$ converges weakly to $Q$ in $L^2(\Omega\times [0,T],\mathcal{S}_2(H))$.

\noindent The fact that $(P,Q)$ is a mild solution of the singular BSDEs \eqref{lbsde} follows  by passing to the limit in equation \eqref{eq:BSDEmildN} .

\medskip
\noindent \textit{Uniqueness:}
Let $(P,Q)$ and $(P',Q')$  two solutions and let
$\bar{P}=P-P'$ and $\bar{Q}=Q-Q'$.
Moreover let $\bar L $ be the weak limit in $L^2_{\mathcal{P}}(\Omega\times [0,T])$, as $N\rightarrow \infty$, of $(T-\cdot)^{\alpha}\sum_{i=1}^NC_i^*(s)\bar Q_s e_i$ and define:
$$\bar{\varphi}_s^N=\bar{L}_s-(T-s)^\alpha\sum_{i=1}^{N} C_i^*(s)\bar{Q}_se_i. $$
Then  equation \eqref{eq:mild} yields:
$$
\bar P_s=-\sum_{i=1}^{\infty}\int_s^T e^{(l-s)A^*}\bar Q_l d\beta_l^i+
 \sum_{i=1}^N \int_t^T  e^{(l-s)A^*} C_i^*(l)\bar Q_le_i dl +\int_s^T  e^{(l-s)A^*}(T-l)^{-\alpha}\bar \varphi^N_l dl.$$
 thus, for any fixed $N$, $(\bar P,\bar Q)$ is the unique mild solution of equation \eqref{eq:exploding} with $\xi_l=(T-l)^{-\alpha}\bar \varphi^N_l$.
By Lemma \ref{lemma:exploding} we obtain that, for all $\mathcal{S}_2$-valued bounded predictable $\Gamma$:
$$\mathbb E\langle \bar{P}_s,x\rangle +\mathbb E\int_s^T \langle \bar{Q}_l,\Gamma_l\rangle_{\mathcal{S}_2(H)}dl=\mathbb E\int_s^T \langle (T-l)^{-\alpha} \tilde{\mathcal{Y}}^{N,\Gamma}_l,\bar{\varphi}_l^N\rangle_{L_2(H)}dl,$$
where again $\tilde{\mathcal{Y}}^{N,\Gamma}$ is the mild solution of equation \eqref{lfsden} with $\gamma=\rho=0$.
By  Theorem \ref{thm:lfsde} we have, in particular:
$\sup_{l\in [t,T], N\in \mathbb{N}} \mathbb E|\tilde{\mathcal{Y}}^{N,\Gamma}_s|^2\le c,$
and $\mathbb E|\tilde{\mathcal{Y}}^{N,\Gamma}_l-{\mathcal{Y}}^{\Gamma}_l|^2\rightarrow 0,\quad \forall l\in [s,T],$
where  ${\mathcal{Y}}^{\Gamma}$ is the mild solution of equation \eqref{lfsden} with $\gamma=\rho=0$. In particular $\tilde{\mathcal{Y}}^{N,\Gamma}$ converges strongly to ${\mathcal{Y}}^{\Gamma}$ in  $L^2_{\mathcal{P}}(\Omega\times [s,T],H)$.

Since $\bar{\varphi}^N$ converges weakly to $0$, we obtain:
$$\mathbb E\langle \bar{P}_s,x\rangle +\mathbb E\int_s^T \langle \bar{Q}_l,\Gamma_l\rangle_{L_2(H)}dl=0,$$
which concludes the proof of uniqueness.  \qed

\subsection{Trace Class Regularity}

We will now prove that the martingale term $Q$ enjoys a \textit{trace class regularity} that will be essential to formulate the maximum principle.

\begin{proposition}\label{prop:trace-class-reg} If $(P,Q)$ is the unique mild solution of equation \eqref{lbsde}, then
$$\mathbb E\int_0^T (T-l)^{2\alpha}|Q_l|^2_{\mathcal{S}_1(H)}dl\le c \left(\mathbb{E}|\eta|^2+\mathbb{E}\int_0^T |f_l|^2 dl \right),
$$
where the constant $c$ only depends on the constants in Hypothesis \ref{Hyp:C_i}.
\end{proposition}

\noindent \textbf{Proof.} Since $Q_l$ is, $\mathbb{P}$-a.s. of class $\mathcal{S}_2(H)$ and therefore compact, it can be written (see Appendix 1) as:
$$Q_l=\sum_{j=1}^\infty a_j(l)h_j(l)\langle g_j(l),\cdot\rangle,$$
where $a_j(l)\in \mathbb R$ and $(h_j(l))_j$, $(g_j(l))_j$ are orthonormal bases in $H$.
Moreover,   we can choose the processes $a_j(l)$, $h_j(l)$,  $ g_j(l)$ ($l\in [0,T]$) to be progressively measurable. Let
$$\Gamma_l^n=\varphi(l) \sum_{j=1}^n sgn(a_j(l))h_j(l)\langle g_j(l),\cdot\rangle ,$$
where $\varphi$ is an arbitrary positive
real-valued bounded progressively measurable process. We note that
$|\Gamma_l^n|_{\mathcal{L}(H)}\le \varphi(l)$ and, being of rank $n$,
the process $\Gamma^n$ is also bounded in $\mathcal{S}_2(H)$.

By \eqref{eq:dualityM} we have
$$\mathbb E\int_0^T \langle Q_l,\Gamma_l^n\rangle_{\mathcal{S}_2(H)}dl
=\mathbb E\langle \mathcal{Y}_T^n,\eta\rangle +\mathbb E\int_0^T \langle \mathcal{Y}^n_l,f_l\rangle dl,$$
where $\mathcal{Y}^n$ is now the solution to equation \eqref{lfsde} with $s=x=\gamma=\rho=0$ and $\Gamma=\Gamma^n$. If we compute $\langle Q_l,\Gamma_l^n\rangle_{\mathcal{S}_2(H)}$
 and we estimate $\mathbb{E}|\mathcal{Y}^n_T|^2$
by  \eqref{eq:gronwall-estimate}  and $\mathbb{E} \int_0^T |\mathcal{Y}^n_l|^2 dl$ by \eqref{eq:gronwall-estimate-int}
we obtain:
\begin{eqnarray*}
  \mathbb E\int_0^T\sum_{j=1}^n |a_j(l)|\varphi(l)d l &\le&
  c\,\left(\mathbb E|\eta|^2\right)^{\frac{1}{2}}\left(\mathbb E\int_0^T(T-l)^{-2\alpha}\varphi(l)^2dl\right)^{\frac{1}{2}}
 \\
    & &  +c\,\mathbb e\left(\int_0^T \mathbb{E}|f_l|^2dl\right)^{\frac{1}{2}}\left(\mathbb E\int_0^T\varphi(l)^2dl\right)^{\frac{1}{2}}.
\end{eqnarray*}
If now $n\rightarrow \infty$, recalling that $\sum_{j=1}^{\infty}|a_j|= |Q|_{\mathcal{S}_1(H)}$
(see Appendix 1 below) we get:
$$ \mathbb{E}\int_0^T \varphi(l)\,|Q_l|_{\mathcal{S}_1}dl\le c_{f,\eta}^{\frac{1}{2}}\left(\mathbb{E}
\int_0^T(T-l)^{-2\alpha}\varphi(l)^2dl\right)^{\frac{1}{2}},$$
where $c_{f,\eta}=c\left(\mathbb{E}|\eta|^2
+\mathbb{E}\int_0^T |f_l|^2 dl \right)$. Denoting
$\tilde{\varphi}(l)=(T-l)^{-\alpha}\varphi(l)$, we can rewrite the last estimate as:
$$\mathbb{E}\int_0^T \tilde{\varphi}(l) [(T-l)^{\alpha} |Q_l|_{L_1(H)}] d l\le  c_{f,\eta}^{\frac{1}{2}}|\tilde{\varphi}|_{L^2_{\mathcal{P}}(\Omega,\times [0,T])}$$
and the claim follows from the arbitrariness of $ \varphi$. \qed

We end this section by proving that, under the following additional assumption, the weak limit that defines the term  $\sum_{i=1}^N \int_s^T e^{(l-s)A^*}C_i^*(l)Q_le_idl$ in Definition \ref{def:solution_BSDE} is indeed a strong limit in $L^1$.
\begin{hypothesis}\label{additional-hyp} We have
 $\sum_{i=1}^\infty |e^{sA^*}C_i^*(l)x|^2\le c s^{-2\alpha}|x|^2$
 for all $ s> 0,$ $l\geq 0$, $x\in H$,  with $\alpha < 1/2$.
\end{hypothesis}

\begin{remark}{\em In the example in Section 2.1 this requirement  coincides with Hypothesis \ref{Hyp:C_i}-{\it 2} since $A$ and  $C_i$ are self adjoint.
}\end{remark}

\begin{proposition}\label{prop:strong} If Hypothesis \ref{additional-hyp} holds, in addition to Hypothesis \ref{Hyp:C_i},
then the sequence $\sum_{i=1}^N  \int_s^T e^{(l-s)A^*}C_i^*(l)Q_le_idl$ converges strongly in $L^1(\Omega,\mathcal{F}_T,\mathbb{P};H)$ as $N\to\infty$. The limit obviously coincides with the weak limit in $L^2(\Omega,\mathcal{F}_T,\mathbb{P};H)$ introduced in Definition \ref{def:solution_BSDE}.
\end{proposition}

\noindent{\bf Proof.}
As above we expand $Q$ as
$Q_l=\sum_{j=1}^\infty a_j(l)h_j(l)\langle g_j(l),\cdot\rangle.$
Then
\begin{eqnarray*}
\sum_{i=1}^N \mathbb E \left| \int_s^T e^{(l-s)A^*}C_i^*(l)Q_le_idl\right |
&\le&\mathbb E  \int_s^T \sum_{i=1}^N \left|e^{(l-s)A^*}C_i^*(l)\sum_{j=1}^\infty a_j(l)h_j(l)\langle g_j(l),e_i\rangle\right|dl\\
&\le& \mathbb E  \int_s^T \sum_{i=1}^N \sum_{j=1}^\infty |a_j(l)| \ |e^{(l-s)A^*}C_i^*(l)h_j(l)| \ |\langle g_j(l),e_i\rangle|dl\\
&\le& \mathbb E  \int_s^T \sum_{j=1}^\infty |a_j(l)|  \left(\sum_{i=1}^\infty |e^{(l-s)A^*}C_i^*(l)h_j(l)|^2\right)^{\frac{1}{2}} \left(\sum_{i=1}^\infty\langle g_j(l),e_i\rangle^2\right)^{\frac{1}{2}}\\
&\le& c \mathbb E  \int_s^T \sum_{j=1}^\infty |a_j(l)|  (l-s)^{-\alpha} dl\\
&=& c \mathbb E  \int_s^T |Q_l|_{\mathcal{S}_1(H)}(l-s)^{-\alpha}dl\\
&\le& c\left( \int_s^T |Q_l|^2_{\mathcal{S}_1(H)}(T-l)^{2\alpha} dl\right)^{\frac{1}{2}}\left(  E  \int_s^T (l-s)^{-2\alpha}(T-l)^{-2\alpha}ds\right)^{\frac{1}{2}}.
\end{eqnarray*}
Strong convergence in $L^1(\Omega;H)$ follows from Proposition \ref{prop:trace-class-reg}. The coincidence  of the two limits is evident by testing them against any bounded $\mathcal{F}_s$-measurable $H$-valued random variable.
\qed

\section{Appendix}

\subsection{Trace Class Operators}\label{appendix:trace-class}
We collect here some basic facts on Hilbert-Schmidt and trace class operators. For a detailed treatment and for the proof of the results stated below see \cite{Ringrose}.

\medskip

\noindent Given a real separable Hilbert space $H$ we denote by {$\mathcal{S}_2(H)$} the Hilbert space of Hilbert Schmidt operators $H\rightarrow H$ endowed with the scalar product
$\<L,M\>_{\mathcal{S}_2(H)}=\sum_{i=1}^{\infty} \<Le_i,Me_i \>_H$ where $(e_i)$ is any orthonormal basis in $H$.

\medskip

\noindent If $L\in \mathcal{S}_2(H)$ then there exists a sequence $(a^L_j)_{j\in \mathbb{N}}\in \ell_2 $ and a pair of orthonormal bases
 $(e^L_j)_{j\in \mathbb{N}}$,  $ (h^L_j)_{j\in \mathbb{N}}$  in $H$ such that
  $ L=\sum_{j=1}^{\infty} a^L_j h^L_j \<e^L_j,\cdot\>$. Moreover $|L|_{\mathcal{S}_2(H)}=\sum_{j=1}^{\infty} (a^L_j)^2.$ Finally if $t\rightarrow L_t$ is a $\mathcal{S}_2(H)$-valued measurable process then the above objects can be selected with the same measurability property.

\medskip

\noindent We define trace class operators in the following way: $ {  \mathcal{S}_1(H)}=\left\{L\in \mathcal{S}_2(H) : |L|_{\mathcal{S}_1}<\infty\right \}$, where
$$\displaystyle|L|_{\mathcal{S}_1}:=\sup\left\{ \<B,L\>_{\mathcal{S}_2}: B \in \mathcal{S}_2(H), \, {|B|_{ \mathcal{L}(H)}\leq 1}\right\}.$$
The following results are true:
\begin{enumerate}
  \item  If $  B\in \mathcal{L}(H)$ and $ L\in \mathcal{S}_1(H)$ then $  LB,\, BL$  are in $   \mathcal{S}_1(H)$, and moreover \newline  $ |LB|_{\mathcal{S}_1(H)}
      \leq |L|_{\mathcal{S}_1(H)} |B|_{\mathcal{L}(H)}$, $ |BL|_{\mathcal{S}_1(H)}\leq |L|_{\mathcal{S}_1(H)} |B|_{\mathcal{L}(H)}$.
  \item If $L\in \mathcal{S}_1(H)$ and $(e_i)$ is an arbitrary orthonormal basis, the trace $\hbox{Tr}(L):= \sum_{i=1}^{\infty}\<e_i, Le_i\>$ converges absolutely and its value is independent of the choice of the basis $(e_i)$.
  \item Using the expansion introduced above $ |L|_{\mathcal{S}_1(H)}=\sum_{j=1}^{\infty} |a^L_j|$, $\hbox{Tr}(L)=\sum_{j=1}^{\infty} a^L_j$ consequently $|\hbox{Tr}(L)|\leq |L|_{\mathcal{S}_1(H)}$.
 \end{enumerate}

\subsection{Linear, Infinite-dimensional BSDEs with Unbounded Terms}
We prove here that infinite dimensional BSDEs can be well-posed even if they include an unbounded term. The following proposition is an extension of the results in \cite{HuPe2} and its proof follows the same lines.
\begin{proposition}\label{appendix:existence_simple_case}
Under the same assumptions and notations of Section 2, let $\eta \in L^2(\Omega, \mathcal{F}_T,\mathbb{P};H)$ and let $\phi$ be a progressively measurable process
in $H$ satisfying
$$\mathbb E\int_0^T (T-s)^{2\alpha}|\varphi_s|^2 ds<\infty,$$
then the BSDE
$$-dp_t=(A^*p_t+\varphi_t)dt-q_tdW_t, \quad p_T=\eta$$
admits a unique mild solution, that is a unique pair of processes $(p,q)$ with $p\in L^2_{\mathcal{P}}(\Omega,{C}([0,T], H))$,
$q\in  L^2_{\mathcal{P}}(\Omega\times[0,T], \mathcal{S}_2(H))$ verifying
$$p_t=e^{(T-t)A^*}\eta+\int_t^T e^{(l-t)A^*}\varphi_l dl-\int_t^T e^{(l-t)A^*}q_l dW_l.$$
Moreover the following estimate holds:
$$\mathbb E\int_t^T |q_s|_{\mathcal{S}_2(H)}^2ds
+\sup_{s\in [t,T]}\mathbb E|p_s|^2\le  c\mathbb{E}|\eta|^2 + c(T-t)^{1-2\alpha}\int_t^T (T-s)^{2\alpha}\mathbb E|\varphi_s|^2ds.$$
\end{proposition}

\noindent{\bf Proof.} The uniqueness  follows directly from \cite{HuPe2} since the difference $(\bar p, \bar q)$ of two solutions satisfies:
$$-d\bar p_t=A^*\bar p_tdt-\bar q_tdW_t, \quad \bar p_T=0.$$
Concerning existence, let us set:
$$p_t= e^{(T-t)A^*}\mathbb E(\eta|{\cal F}_t)+\int_t^T  e^{(s-t)A^*}\mathbb E(\varphi_s|{\cal F}_t)ds.$$
Moreover, by the martingale representation theorem,
$$\mathbb E(\varphi_s|{\cal F}_t)=\varphi_s-\int_t^s g(s,l)dW_l, \quad \mathbb E(\eta|{\cal F}_t)=\eta-\int_t^T h(l)dW_l.$$
Notice that
$$\mathbb E\int_t^\rho |g(\rho,\sigma)|^2d\sigma\le \mathbb E|\varphi_\rho|^2+\mathbb E|\mathbb E(\varphi_\rho|{\cal F}_t)|^2\le 2\mathbb E|\varphi_\rho|^2,\quad \mathbb E\int_t^\rho |h(\sigma)|^2d\sigma\le 2\mathbb E|\eta|^2.$$
From the above two equations, we have:
\begin{eqnarray*}
p_t&=&\int_t^T e^{(s-t)A^*}\varphi_sds-\int_t^T  e^{(s-t)A^*}\left(\int_t^s g(s,l)dW_l\right)ds + e^{(T-t)A^*}\eta- \int_t^T  e^{(T-t)A^*}h(l)dW_l\\
&=&\int_t^T e^{(s-t)A^*}\varphi_sds-\int_t^T  e^{(l-t)A^*}\left(\int_l^T e^{(s-l)A^*} g(s,l)ds \right)dW_l\\
&&\quad + e^{(T-t)A^*}\eta- \int_t^T  e^{(l-t)A^*}e^{(T-l)A^*}h(l)dW_l.\end{eqnarray*}

So setting
$$q_l=\int_l^T e^{(s-l)A^*} g(s,l)ds +e^{(T-l)A^*}h(l),$$
we deduce that $(P,Q)$ is the unique solution

Let us now establish the estimates.
\begin{eqnarray*}\mathbb E|q_\sigma|^2&\leq& c\left(\int_\sigma^T (T-\rho)^{-\alpha}(T-\rho)^{\alpha}\mathbb E|g(\rho,\sigma)|d\rho\right)^2+c\mathbb{E}|h(\sigma)|^2 d \sigma\\&\le& (T-t)^{1-2\alpha}\int_\sigma^T (T-\rho)^{2\alpha}\mathbb E|g(\rho,\sigma)|^2d\rho + c\mathbb{E}|h(\sigma)|^2 d\sigma.\end{eqnarray*}

Thus
\begin{eqnarray*}
\int_t^T \mathbb E|q_\sigma|^2d\sigma&\le& (T-t)^{1-2\alpha}\int_t^T \int_\sigma^T (T-\rho)^{2\alpha}\mathbb E|g(\rho,\sigma)|^2d\rho d\sigma+c\int_t^T \mathbb{E}|h(\sigma)|^2 d\sigma \\
&\le& (T-t)^{1-2\alpha}\int_t^T (T-\rho)^{2\alpha}\left( \int_t^\rho E|g(\rho,\sigma)|^2d\sigma\right)d\rho+c\int_t^T \mathbb{E}|h(\sigma)|^2 d\sigma \\
&\le& (T-t)^{1-2\alpha}\int_t^T (T-\rho)^{2\alpha}\mathbb E|\varphi(\rho)|^2 d\rho+c \mathbb{E}|\eta|^2 .
\end{eqnarray*}

On the other hand,
\begin{eqnarray*}
\mathbb E|p_t|^2&\leq& c\mathbb E\left|\int_t^T \mathbb E(\varphi_s|{\cal F}_t)ds\right|^2+\mathbb E| \mathbb E(\eta|{\cal F}_t)|^2
=\mathbb E\left[\int_t^T (T-s)^{-\alpha}(T-s)^{\alpha}|\mathbb E (\varphi_s|{\cal F}_t)|ds\right]^2+\mathbb{E}|\eta|^2\\
&\le& (T-t)^{1-2\alpha}\int_t^T (T-s)^{2\alpha}\mathbb E |\varphi(s)|^2 ds +\mathbb{E}|\eta|^2,
\end{eqnarray*}
which concludes the proof.
\qed

\noindent{\bf Proof of Lemma \ref{lemma:exploding}.} Using Proposition \ref{appendix:existence_simple_case}, existence and uniqueness  follows by a standard contraction argument (see \cite{HuPe2}).
The final duality property is  established by a simple truncation argument.
\qed


\begin{thebibliography}{11}


\bibitem{Be}
A. Bensoussan.
 Stochastic maximum principle for distributed parameter systems.
J. Franklin Inst. 315 (1983), no. 5-6, 387--406.



\bibitem{DaZa} G. Da Prato, J.  Zabczyk.
Stochastic equations in infinite dimensions.
Encyclopedia of Mathematics and its Applications, 44. Cambridge University Press, Cambridge,
 1992.

 \bibitem{DaZa2} G. Da Prato, J.  Zabczyk. Ergodicity for Infinite Dimensional Systems.
 London Mathematical Society Lecture Note Series, 229.
 Cambridge University Press, Cambridge, 1996.


\bibitem{DuMe}
K. Du, Q. Meng.
Stochastic maximum principle for infinite dimensional control systems. Preprint
  arXiv:1208.0529, 2012.

\bibitem{DuMe2}
K. Du, Q. Meng.
A maximum principle for optimal control of stochastic evolution equations. SIAM J. Control Optim. 51
(2013), no. 6, 4343--4362.

 \bibitem{FuHuTe}
 M. Fuhrman, Y.  Hu, G. Tessitore.
 Stochastic maximum principle for optimal control of SPDEs.
 Appl. Math. Optim. 68 (2013), no. 2, 181--217.

\bibitem{FT1} M. Fuhrman, G.  Tessitore. Nonlinear Kolmogorov equations in infinite dimensional spaces: the backward stochastic differential equations approach and applications to optimal control. Ann. Probab. 30 (2002), no. 3, 1397--1465.

\bibitem{Gu}
G. Guatteri.
 Stochastic maximum principle for SPDEs with noise and control on the boundary.
 Systems Control Lett. 60 (2011), no. 3, 198--204.

 \bibitem{Gu-Te}
G. Guatteri, G. Tessitore.
On the backward stochastic Riccati equation in infinite
               dimensions. SIAM J. Control and Optim. 
44 (2005), no. 1, 159--194.

\bibitem{Henry} D. Henry.
Geometric theory of semilinear parabolic equations.
Lecture Notes in Mathematics, 840. Springer, Berlin, 1981.


\bibitem{HuPe} Y. Hu, S. Peng.
Adapted solution of a backward semilinear stochastic evolution equation.
Stochastic Anal. Appl. 9 (1991), no. 4, 445--459.

\bibitem{HuPe2}
Y. Hu, S. Peng.
Maximum principle for semilinear stochastic evolution control systems.
Stochastics Stochastics Rep. 33 (1990), no. 3-4, 159--180.

\bibitem{lu} A. Lunardi. Analytic semigroups and optimal
regularity in parabolic problems.
Progress in Nonlinear Differential Equations and their Applications, 16.
Birkhauser, Basel, 1995.


\bibitem{LuZh}
Q. L\"u, X. Zhang.
General Pontryagin-type stochastic maximum principle and
backward stochastic evolution equations in infinite dimensions. 
SpringerBriefs in Mathematices. Springer, Cham, 2014.


\bibitem{OkSuZh}
B. \O ksendal, A.  Sulem, T.  Zhang.
 Optimal control of stochastic delay equations and time-advanced backward
 stochastic differential equations.
 Adv. in Appl. Probab. 43 (2011), no. 2, 572--596.


\bibitem{pa}  A. Pazy. Semigroups of linear operators and
applications to partial differential equations.
Applied Mathematical Sciences, 44. Springer, New York, 1983.


\bibitem{Pe} S. Peng.
 A general stochastic maximum principle for optimal control problems.
SIAM J. Control Optim. 28 (1990), no. 4, 966--979.
\bibitem{Ringrose}

J. R. Ringrose.
Compact non-self-adjoint operators, Van Nostrand, London, 1971.


\bibitem{TaLi}
S. Tang, X. Li.
 Maximum principle for optimal control of distributed parameter
 stochastic systems with random jumps.
 Differential equations, dynamical systems, and control science, 867--890,
 Lecture Notes in Pure and Appl. Math., 152, Dekker, New York,  1994.

\bibitem{Te1} G. Tessitore. Some remarks on the Riccati equation arising in an optimal control problem with state- and control-dependent noise.  SIAM J. Control Optim. 30 (1992), no. 3, 717--744.

\bibitem{Te2} G. Tessitore.
Existence, uniqueness and space regularity of the adapted solutions of a backward SPDE.
Stochastic Anal. Appl. 14 (1996), no 4, 461--486.



 \bibitem{vNVW}  J. M. A. M. van Neerven, M. C. Veraar, L. Weis. Stochastic integration in UMD Banach spaces. Ann. Probab. 35  (2007), no. 4,  1438--1478.

\bibitem{YoZh}  J. Yong, X.Y.  Zhou.
Stochastic controls.
Hamiltonian systems and HJB equations.
Applications of Mathematics (New York), 43. Springer, New York, 1999.

\bibitem{Zh}
X.Y.  Zhou.
On the necessary conditions of optimal controls for stochastic
partial differential equations. SIAM J. Control Optim. 31 (1993), no. 6, 1462--1478.
\end{thebibliography}
\end{document}